\documentclass[11pt]{article}
\usepackage[utf8]{inputenc}
\usepackage{graphicx} 
\usepackage{amsmath, amsthm, amssymb} 
\usepackage{hyperref} 
\usepackage{multicol}
\usepackage{multirow}
\usepackage{caption}
\usepackage{subcaption}
\usepackage[labelfont=bf]{caption}
\usepackage{pdflscape}
\usepackage{xcolor}

\usepackage{algorithm}
\usepackage[noend]{algpseudocode}

\newcommand{\myalgsheader}[0]

\algnewcommand{\IIf}[1]{\State\algorithmicif\ #1\ \algorithmicthen}
\algnewcommand{\EndIIf}{\unskip\ \algorithmicend\ \algorithmicif}
\algnewcommand{\IElse}[1]{\State\algorithmicelse\ #1\ }
\algnewcommand{\IfThenElse}[3]{
  \State \algorithmicif\ #1\ \algorithmicthen\ #2\ \algorithmicelse\ #3}

\begin{document}
\title{Multistage Mixed Precision Iterative Refinement\footnote{We acknowledge funding from the Charles University PRIMUS project
No. PRIMUS/19/SCI/11, Charles University Research Program No. UNCE/SCI/023, and the Exascale Computing Project (17-SC-20-SC), a collaborative effort of the U.S. Department of Energy Office of Science and the National Nuclear Security Administration.}}

\author{Eda Oktay and Erin Carson}
\date{}
\maketitle

\paragraph{Abstract.}
Low precision arithmetic, in particular half precision (16-bit) floating point arithmetic, is now available in commercial hardware. 
Using lower precision can offer significant savings in computation and communication costs with proportional savings in energy. Motivated by this, there has been a renewed interest in mixed precision iterative refinement schemes for solving linear systems $Ax=b$, and new variants of GMRES-based iterative refinement have been developed. Each particular variant with a given combination of precisions leads to different condition number-based constraints for convergence of the backward and forward errors, and each has different performance costs. The constraints for convergence given in the literature are, as an artifact of the analyses, often overly strict in practice, and thus could lead a user to select a more expensive variant when a less expensive one would have sufficed.

In this work, we develop a multistage mixed precision iterative refinement solver which aims to combine existing mixed precision approaches to balance performance and accuracy and improve usability. For a user-specified initial combination of precisions, the algorithm begins with the least expensive approach and convergence is monitored via inexpensive computations with quantities produced during the iteration. If slow convergence or divergence is detected using particular stopping criteria, the algorithm switches to use a more expensive, but more reliable variant. A novel aspect of our approach is that, unlike existing implementations, our algorithm first attempts to use ``stronger'' GMRES-based solvers for the solution update before resorting to increasing the precision(s). In some scenarios, this can avoid the need to refactorize the matrix in higher precision. We perform extensive numerical experiments on a variety of random dense problems and problems from real applications which confirm the benefits of the multistage approach.

\section{Introduction} \label{sec:first_intro}
Iterative refinement (IR) is frequently used in solving linear systems $Ax = b$, where $A \in \mathbb{R}^{n \times n}, x,b \in \mathbb{R}^n$, to improve the accuracy of a
computed approximate solution $\hat{x} \in \mathbb{R}^n$. Typically, one computes an initial approximate solution $\hat{x}_0\in \mathbb{R}^n$ using Gaussian
elimination with partial pivoting (GEPP), saving the approximate factorization $A \approx \hat{L}\hat{U}$, where $\hat{L} \in \mathbb{R}^{n \times n}$ is a lower triangular matrix with unit diagonal, and $\hat{U} \in \mathbb{R}^{n \times n}$ is an upper triangular matrix. After computing the
residual $\hat{r}=b-A\hat{x} \in \mathbb{R}^n$ (potentially in higher precision), one reuses $\hat{L}$ and $\hat{U}$ to solve the system $A \hat{d} = \hat{r}$, where $\hat{d} \in \mathbb{R}^n$. The
original approximate solution is subsequently refined by adding the corrective term, $\hat{x} = \hat{x} + \hat{d}$. This process
can be repeated until either the desired accuracy is reached or the iterative refinement process converges to an
approximate solution.
Iterative refinement algorithms that exploit multiple different precisions have seen renewed attention recently due to the emergence of mixed precision hardware, and form the basis for the recently developed HPL-AI benchmark, on which today's top supercomputers exceed exascale performance; see e.g., \cite{hplai,k:20,top500}.

Algorithm \ref{alg:sir} shows a general mixed precision iterative refinement scheme. Following the work of \cite{ch:18}, there are four different precisions specified: $u_f$ denotes the factorization precision, $u$ denotes the working precision, $u_r$ denotes the precision for the residual computation, and $u_s$ denotes the effective precision with which the correction equation is solved. This final precision, $u_s$, is not a hardware precision but will rather depend on the particular solver used in line \ref{corrsolve} and the particular precision(s) used within the solver (which may be different from $u$, $u_f$, and $u_r$). It is assumed that $u_f\geq u\geq u_r$. 

\begin{algorithm}[htbp!]
	\caption{General Iterative Refinement Scheme \label{alg:sir}} 
	\begin{algorithmic}[1]
		\Require{$n \times n$ matrix $A$; right-hand side $b$; maximum number of refinement steps $i_{max} \in \mathbb{N}^+$.}
		\Ensure{Approximate solution $x_{i+1}$ to $Ax = b$.}
		\State{Compute LU factorization $A = LU$ in precision $u_f$}.
		\State{Solve $Ax_0 = b$ by substitution in precision $u_f$; store $x_0$ in precision $u$.}
		\For{$i = 0$: $i_{max} - 1$}    
			\State{Compute $r_i = b - Ax_i$ in precision $u_r$; store in precision $u$.}
			\State{Solve $Ad_{i+1} = r_i$ in precision $u_s$; store $d_{i+1}$ in precision $u$. \label{corrsolve}}
			\State{Compute $x_{i+1} = x_i + d_{i+1}$ in precision $u$.}
			\IIf{converged} return $x_{i+1}$. \EndIIf
		\EndFor
	\end{algorithmic}
\end{algorithm}

The choice of precisions $u$, $u_f$, and $u_r$, and the choice of a solver to use in line \ref{corrsolve} defines different variants of iterative refinement. We refer to any variant that solves the correction equation in line \ref{corrsolve} using triangular solves with the computed LU factors as ``standard iterative refinement'' (SIR). For SIR, the effective precision of the correction solve is limited by the precision with which the factorization is computed, and thus we have $u_s = u_f$. 
The most commonly-used SIR variant is what we call ``traditional'' iterative refinement, in which $u_f=u$ and $u_r=u^2$. For instance, if the working precision is single, i.e., $u = 5.96\cdot 10^{-8}$, then for $u_r=u^2 \approx 10^{-16}$, double precision is used. Traditional iterative refinement was used already by Wilkinson in 1948 \cite{w:48}, and was analyzed in fixed point arithmetic by Wilkinson in 1963 \cite{w:63} and in floating point arithmetic by Moler in 1967 \cite{m:67}. There have also been analyses of fixed precision iterative refinement, in which $u_f=u=u_r$ \cite{jw:77}, as well as low precision factorization iterative refinement, in which $u_f^2=u=u_r$ \cite{lllkbd:06}. Motivated by the trend of low and mixed precision capabilities in hardware, the authors in \cite{ch:18} developed and analyzed a three-precision iterative refinement scheme, in which $u_f$, $u$, and $u_r$ may differ, which generalizes (and in some cases, improves the bounds for) many existing variants. 
For references to analyses of variants of iterative refinement, see \cite[Table 1.1]{ch:18}. 

If $A$ is very ill conditioned or badly scaled, SIR can fail, i.e., the error is not eventually bounded by a small multiple of machine precision. In extreme cases of ill conditioning, the error can grow with each refinement step. In \cite{ch:17} and \cite{ch:18}, the authors developed a mixed precision GMRES-based iterative refinement scheme (GMRES-IR), shown in Algorithm \ref{alg:gmresir}. The only difference with SIR is the way in which the correction equation is solved. The idea is that instead of using the LU factors to solve for the correction in each step, one can instead use these factors as (left) preconditioners for a preconditioned GMRES method which solves for the correction. In this way, the effective solve precision becomes $u_s=u$, and more ill-conditioned problems can be handled relative to SIR. The GMRES-based refinement approaches have seen much success in practice. Current GMRES-based iterative refinement variants are implemented in the MAGMA library (2.5.0) and in the NVIDIA cuSOLVER library. GMRES-IR also forms the basis for the new HPL-AI benchmark, used to rank computers in the TOP500 list \cite{hplai}. Experiments detailing the performance benefits of three-precision GMRES-IR and SIR can be found in \cite{ht:18}. 

\begin{algorithm}[htbp!]
	\caption{GMRES-IR \cite{ch:17} \label{alg:gmresir}}
	\begin{algorithmic}[1]
		\Require{$n \times n$ matrix $A$; right-hand side $b$; maximum number of refinement steps $i_{max}$; GMRES convergence tolerance $\tau$.}
		\Ensure{Approximate solution $\hat{x}$ to $Ax = b$.}
		\State{Compute LU factorization $A = LU$ in precision $u_f$.}
		\State{Solve $Ax_0 = b$ by substitution in precision $u_f$; store $x_0$ in precision $u$.}
		\For{$i = 0$: $i_{max} - 1$}   
			\State{Compute $r_i = b - Ax_i$ in precision $u_r$; store in precision $u$.}
			\State{Solve $U^{-1}L^{-1}Ad_{i+1}=U^{-1}L^{-1}r_i$ by GMRES in working precision $u$, with matrix-vector products with $\tilde{A}=U^{-1}L^{-1}A$ computed at precision $u^2$; store $d_{i+1}$ in precision $u$. \label{corrsolveg}}
			\State{Compute $x_{i+1}=x_i+d_{i+1}$ in precision $u$.}
			\IIf{converged} return $x_{i+1}$. \EndIIf
		\EndFor
	\end{algorithmic}
\end{algorithm}

Although GMRES-IR can succeed where SIR fails, each step of GMRES-IR is potentially more expensive
than each step of SIR (albeit still potentially less expensive than running the entire process in
double the working precision). Whereas each SIR refinement step involves only two triangular solves in precision $u$, each GMRES-IR step involves two triangular solves in precision $u^2$ in \emph{each} GMRES iteration (in addition to a matrix-vector product in precision $u^2$ as well as multiple vector operations in precision $u$). The convergence behavior of each GMRES solve thus plays a large role, and it is important for the performance of GMRES-IR that each call to GMRES converges relatively quickly. It is additionally a performance concern that each GMRES iteration requires computations in higher precision; in order to obtain the needed bound on backward error of the GMRES solve, the authors in \cite{ch:17, ch:18} required that the preconditioned system $U^{-1}L^{-1}A$ (not formed explicitly) is applied to a vector in precision $u^2$. 

The authors of \cite{ch:17,ch:18} suggested that this required use of extra precision was likely to be overly strict in many practical scenarios. Indeed, most practical implementations of GMRES-based iterative refinement do not use this extra precision \cite{h:21}. This motivated Amestoy et al. \cite{h:21} to develop an extension of the GMRES-IR algorithm, called GMRES-IR5. The authors in \cite{h:21} revisit the proof of backward stability for GMRES from \cite{prs:06} and develop a bound on the backward error for a mixed precision GMRES algorithm, where $u_g \geq u$ is the working precision used within GMRES and $u_p$ is the precision in which the preconditioned matrix is applied to a vector. This enables a five-precision GMRES-IR scheme, where $u, u_f, u_r, u_g$, and $u_p$ may take on different values. A particular improvement over the GMRES-IR scheme of \cite{ch:17, ch:18} is that the analysis provides  bounds on forward and backward errors for a variant of GMRES-IR in which the entire GMRES iteration is carried out in a uniform precision, i.e., the analysis does not require that extra precision is used in applying the preconditioned matrix to a vector. This particular variant of the algorithm is thus less expensive than the former GMRES-IR in terms of both time and memory. The drawback is that it is only theoretically applicable to a smaller set of problems due to a tighter constraint on condition number. We call the particular instance of GMRES-IR5 where $u=u_g=u_p$ ``SGMRES-IR'' (for ``simpler''), shown in Algorithm \ref{alg:sgmresir}. To make more precise the relative costs of each step of SIR, SGMRES-IR, and GMRES-IR, we list their costs in terms of asymptotic computational complexity in Table \ref{tab:cost}. We discuss the constraints on condition number under which each variant converges later in Section \ref{sec:errorbounds}.

\begin{algorithm}[htbp!]
	\caption{SGMRES-IR (a particular variant of GMRES-IR5 \cite{h:21}) \label{alg:sgmresir}}
	\begin{algorithmic}[1]
		\Require{$n \times n$ matrix $A$; right-hand side $b$; maximum number of refinement steps $i_{max}$; GMRES convergence tolerance $\tau$.}
		\Ensure{Approximate solution $x_{i+1}$ to $Ax = b$.}
		\State{Compute LU factorization $A = LU$ in precision $u_f$.}
		\State{Solve $Ax_0 = b$ by substitution in precision $u_f$; store $x_0$ in precision $u$.}
		\For{$i = 0$: $i_{max} - 1$}   
			\State{Compute $r_i = b - Ax_i$	in precision $u_r$; store in precision $u$.}
			\State{Solve $U^{-1}L^{-1}Ad_{i+1}=U^{-1}L^{-1}r_i$ by GMRES in working precision $u$, with matrix-vector products with $\tilde{A}=U^{-1}L^{-1}A$ computed at precision $u$; store $d_{i+1}$ in precision $u$. \label{corrsolvesg}}
			\State{Compute $x_{i+1}=x_i+d_{i+1}$ in precision $u$.}
			\IIf{converged} return $x_{i+1}$. \EndIIf
		\EndFor
	\end{algorithmic}
\end{algorithm}

\begin{table}[]
\centering
\caption{Asymptotic computational complexity of operations in each refinement step for SIR, SGMRES-IR, and GMRES-IR.}
\label{tab:cost}
\begin{tabular}{|c|lll|}
\hline
\begin{tabular}[c]{@{}c@{}}Once per IR solve \\ (all variants)\end{tabular} & $O(n^3)$ & in precision $u_f$ & (LU fact.) \\ \hline
SIR step & $O(n^2)$ & in precision $u_f$ & (tri. solves) \\ \hline
\multirow{3}{*}{\begin{tabular}[c]{@{}c@{}}SGMRES-IR step\\ ($k$ GMRES iterations)\end{tabular}} & $O(nk^2)$ & in precision $u$ & (orthog.) \\
 & $O(nnz\cdot k)$ & in precision $u$ & (SpMV) \\
 & $O(n^2k)$ & in precision $u$ & (precond.) \\ \hline
\multirow{3}{*}{\begin{tabular}[c]{@{}c@{}}GMRES-IR step\\ ($k$ GMRES iterations)\end{tabular}} & $O(nk^2)$ & in precision $u$ & (orthog.) \\
 & $O(nnz\cdot k)$ & in precision $u^2$ & (SpMV) \\
 & $O(n^2k)$ & in precision $u^2$ & (precond.) \\ \hline
\multirow{2}{*}{\begin{tabular}[c]{@{}c@{}}Once per refinement step\\ (all variants)\end{tabular}} & $O(nnz)$ & in precision $u_r$ & (residual comp.) \\
 & $O(n)$ & in precision $u$ & (sol. update) \\ \hline
\end{tabular}
\end{table}

In terms of both cost and range of condition numbers to which the algorithm can be applied, we expect SGMRES-IR to be, in general, somewhere between SIR and GMRES-IR. For example, if we use single precision for $u_f$, double precision for $u$, and quadruple precision for $u_r$, SIR is guaranteed to converge to the level of double precision in both forward and backward errors as long as the infinity-norm condition number of the matrix $A$, $\kappa_\infty(A)=\Vert A \Vert_{\infty} \Vert A^{-1} \Vert_{\infty}$, is less than $2 \cdot 10^7$. For SGMRES-IR, this constraint becomes $\kappa_\infty(A)\leq 10^{10}$. GMRES-IR, on the other hand, only requires $\kappa_\infty(A)\leq 2\cdot 10^{15}$. 
Thus going from SIR to SGMRES-IR to GMRES-IR, we expect that these algorithms will be increasingly expensive but also expect that they will converge for increasingly ill-conditioned matrices. We note this may not always be the case. For some precision combinations, SGMRES-IR has a tighter constraint on condition number than SIR (see \cite{h:21}), although this may be an artifact of the analysis. Also, the metric of relative ``cost'' is difficult to determine a priori, in particular between SGMRES-IR and GMRES-IR, since it depends on the number of GMRES iterations required in each refinement step. 

It is thus difficult to choose a priori which particular variant of iterative refinement is the most appropriate for a particular problem. Even if the matrix condition number meets the constraint for convergence for the chosen iterative refinement algorithm, the convergence rate may be unacceptably slow, or each step may require so many GMRES iterations that it becomes impractical. Further, as our experiments show, the condition number constraints in the literature are can be too tight, meaning that for some problems a given refinement scheme converges even when the analysis indicates that it may not. This could lead users to select a more expensive iterative refinement variant than is actually needed in practice.

In this work, we aim to solve this problem through the development of a multistage, three-precision iterative refinement scheme, which we abbreviate MSIR. Our approach automatically switches between solvers and precisions if slow convergence (of the refinement scheme itself or of the inner GMRES solves) is detected using stopping criteria adapted from the work in \cite{dh:06}. Two novel aspects of our approach are (1) we attempt to use ``stronger'' solvers before resorting to increasing the precision of the factorization, and (2) when executing a GMRES-based refinement algorithm, we modify the stopping criteria to also restrict the number of GMRES iterations per refinement step. 

Table \ref{tab:cost} shows why first switching the solver may be more favorable from a performance perspective; whereas increasing the precision and recomputing the factorization will cost $O(n^3)$ flops in precision $u_f^2$, where $u_f$ is the current factorization precision, using the existing factorization and performing $k$ total GMRES iterations may be faster, requiring $O(n^2k)$ flops in precision $u$. This motivates point (2) above. If the number of GMRES iterations $k$ is too large, then one (S)GMRES-IR refinement step at least as expensive as recomputing the factorization in a higher precision ($u$ for SGMRES-IR or $u^2$ for GMRES-IR).
 
Our approach may serve to improve existing multistage iterative refinement implementations. For example, the MAGMA library \cite{tdb10} currently uses a variant of SGMRES-IR and if convergence is not detected after the specified maximum number of iterations, the factorization is recomputed in a higher precision and the refinement restarts. Our numerical experiments confirm that it may be beneficial to first try a different solver before resorting to recomputing the factorization.

In Section \ref{sec:tsir}, we present the MSIR algorithm and give a motivating numerical example. We then give details of the stopping criteria used and summarize the analysis for each algorithm variant from \cite{ch:18} and \cite{h:21}. In Section \ref{sec:results} we present more thorough numerical experiments on both random dense matrices and matrices from the SuiteSparse collection \cite{dh:11}. We conclude and discuss future extensions in Section \ref{sec:conclusions}.




\section{The MSIR algorithm} 
\label{sec:tsir}
In order to balance reliability and cost, we develop a multistage iterative refinement algorithm (MSIR), presented in Algorithm \ref{alg:tsir}. The algorithm starts with three-precision SIR (as in \cite{ch:18}) and, using the stopping criteria developed in \cite{dh:06}, switches to SGMRES-IR if the algorithm is not converging at an acceptable rate (or not converging at all). Then, using the same stopping criteria along with an additional constraint on the number of GMRES iterations per refinement step, the algorithm may choose to switch a second time to the GMRES-IR algorithm (which uses higher precision in applying the preconditioned matrix to a vector within GMRES). If convergence is still not achieved or is too slow, then as a fail-safe, we increase the factorization precision $u_f$ (and the other precisions if necessary to satisfy $u_f\geq u$, $u_r \leq u^2$), recompute the LU factorization, and begin the process again with SIR. We enforce $u_r \leq u^2$ in order to guarantee the convergence of the forward error to the level of the working precision, but this strategy for increasing precisions could be modified in practice. For instance, if the initial setting is $(u_f,u,u_r)$ = (half, single, double), we would double only $u_f$ and use $(u_f,u,u_r)$ = (single, single, double), and then if we need to increase precisions again, we would use $(u_f,u,u_r)$ = (double, double, quad). 


\begin{algorithm}[htbp!]
	\caption{Multistage Iterative Refinement (MSIR) \label{alg:tsir}}
	\scriptsize
	\begin{algorithmic}[1]
	\Require{$n\times n$ matrix $A$; right-hand-side $b$; maximum number of refinement steps of each type $i_{max}$; GMRES convergence tolerance	$\tau$; stopping criteria parameter $\rho_{thresh}$; maximum GMRES iterations $k_{max} \in \mathbb{N}^+$; initial factorization precision $u_f$; initial working precision $u$; initial residual precision $u_r$.}
	\Ensure{Approximate solution $x_{i+1}$ to $Ax = b$, boolean \texttt{cged}.}
	\State{Compute LU factorization $A = LU$ in precision $u_f$.} \label{tsir:lufact}
	\State{Solve $Ax_0 = b$ by substitution in precision $u_f$; store $x_0$ in precision $u$. \label{tsir:initsolve}}
	\State{Initialize: $d_{0} = \infty$; \texttt{alg} $=$ \texttt{SIR}; \texttt{iter} $=0$; $i=0$; \texttt{cged} $=0$; $\rho_{max} = 0$.}
	\While{not \texttt{cged}}  
		\State{Compute $r_i = b - Ax_i$ in precision $u_r$; Scale $r_i = r_i/\|r_i\|_\infty$; store in precision $u$. \label{r1}}
		\If{\texttt{alg} $=$ \texttt{SIR}} 
			\State{\texttt{iter} $=$ \texttt{iter} $+1$}
			\State{Compute $d_{i+1}=U^{-1}(L^{-1}r_i)$ in precision $u_f$; store $d_{i+1}$ in precision $u$.}
			\IIf{$d_{i+1}$ contains \texttt{Inf} or \texttt{NaN}} 
					\texttt{alg} $=$ \texttt{SGMRES-IR}; \texttt{iter} $=0$; break. \EndIIf \label{tsir:nancheck}
			\State{Compute $x_{i+1}=x_i+\|r_i\|_\infty d_{i+1}$ in precision $u$. \label{x1}}
			\State{$z=\|d_{i+1}\|_\infty/\|x_i\|_\infty$; \quad $v=\|d_{i+1}\|_\infty/\|d_{i}\|_\infty$; \quad $\rho_{max}= \max(\rho_{max},v)$; \quad $\phi_i = z/(1-\rho_{max})$ }
				\If{$z\leq u$ or $v\geq \rho_{thresh}$ or \texttt{iter} $> i_{max}$ or $\phi_i \leq \sqrt{n}u$}
						\If {not converged} \label{conv1}
								\State{\texttt{alg} $=$ \texttt{SGMRES-IR}; \texttt{iter} $=0$.} 
								\IIf {$\phi_{i} > \phi_0$} $x_{i+1} = x_0$ \EndIIf \label{resetx}
					   \Else
								\State{\texttt{cged} $=1$} 
						 \EndIf
				\EndIf
		\ElsIf{\texttt{alg} $=$ \texttt{SGMRES-IR}}
			\State{\texttt{iter} $=$ \texttt{iter} $+1$}
			\State{Solve $U^{-1}L^{-1}Ad_{i+1}=U^{-1}L^{-1}r_i$ by GMRES in precision $u$ with matrix-vector products with $\tilde{A}=U^{-1}L^{-1}A$ computed at precision $u$; store $d_{i+1}$ in precision $u$.}
			\State{Compute $x_{i+1}=x_i+\|r_i\|_\infty d_{i+1}$ in precision $u$. \label{x2}}
			\State{$z=\|d_{i+1}\|_\infty/\|x_i\|_\infty$; \quad $v=\|d_{i+1}\|_\infty/\|d_i\|_\infty$; \quad $\rho_{max}= \max(\rho_{max},v)$; \quad $\phi_i = z/(1-\rho_{max})$}
			\If{$z\leq u$ or $v\geq \rho_{thresh}$ or $\texttt{iter} > i_{max}$ or $k_{GMRES} > k_{max}$ or $\phi_i \leq \sqrt{n}u$}
						\If {not converged} \label{conv2}
								\State{\texttt{alg} $=$ \texttt{GMRES-IR}; \texttt{iter} $=0$.} 
								\IIf {$\phi_{i} > \phi_0$} $x_{i+1} = x_0$ \EndIIf \label{resetx2}
					   \Else
								\State{\texttt{cged} $=1$} 
						 \EndIf
			\EndIf	
		\ElsIf{\texttt{alg} $=$ \texttt{GMRES-IR}}
			\State{\texttt{iter} $=$ \texttt{iter} $+1$}
			\State{Solve $U^{-1}L^{-1}Ad_{i+1}=U^{-1}L^{-1}r_i$ by GMRES in precision $u$ with matrix-vector products with $	\tilde{A}=U^{-1}L^{-1}A$ computed at precision $u^2$; store $d_{i+1}$ in precision $u$.}
			\State{Compute $x_{i+1}=x_i+\|r_i\|_\infty d_{i+1}$ in precision $u$. \label{x3}}
			\State{$z=\|d_{i+1}\|_\infty/\|x_i\|_\infty$; \quad $v=\|d_{i+1}\|_\infty/\|d_i\|_\infty$; \quad $\rho_{max}= \max(\rho_{max},v)$; \quad $\phi_i = z/(1-\rho_{max})$}
			\If{$z\leq u$ or $v\geq \rho_{thresh}$ or $\texttt{iter} > i_{max}$ or $k_{GMRES} > k_{max}$ or $\phi_i \leq \sqrt{n}u$}
			\If {not converged} \label{conv3}
			\State{$u_f = u_f^2$; \texttt{alg} $=$ \texttt{SIR}; \texttt{iter} $=0$.} 
			\State{Compute LU factorization $A = LU$ in precision $u_f$.}
			\IIf{$u_f< u$} $u = u_f$ \EndIIf
			\IIf{$u_r> u^2$} $u_r = u^2$ \EndIIf
			\IIf {$\phi_{i} > \phi_0$} $x_{i+1} = x_0$ \EndIIf \label{resetx3}
			\Else
			\State{\texttt{cged} $=1$} 
			\EndIf
			\EndIf
		\EndIf	
		\State{$i = i+1$}
	\EndWhile
\end{algorithmic}
\end{algorithm}

Before explaining the details of the algorithm, we begin with a brief motivating example illustrating how MSIR works. We restrict ourselves to IEEE precisions and use initial precisions $(u_f,u,u_r)$ = (single, double, quad) and test random dense matrices of size $100 \times 100$ generated using the MATLAB command \verb|gallery('randsvd',n,kappa(i),2)|, where \verb|kappa| is the array of the desired 2-norm condition numbers. Here we test 2-norm condition numbers $10^1, 10^5$, and $10^{16}$. We set $b$ to be a vector of normally distributed numbers generated by the MATLAB command \verb|randn|. For reproducibility, we use the MATLAB command \verb|rng(1)| to seed the random number generator before generating each linear system.  The GMRES convergence tolerance $\tau \in \mathbb{R}^+$, which also appears in Algorithms \ref{alg:gmresir} and \ref{alg:sgmresir}, dictates the stopping criterion for the inner GMRES iterations. The algorithm is considered to be converged if the relative (preconditioned) residual norm drops below $\tau$. In all tests here we use $\tau=10^{-6}$. The particular criteria by which we switch between solvers is discussed in detail in Section \ref{convcriteria}.

\begin{figure}
\centering
	\includegraphics[trim={3cm 8cm 4cm 8cm},clip, width=.45\textwidth]{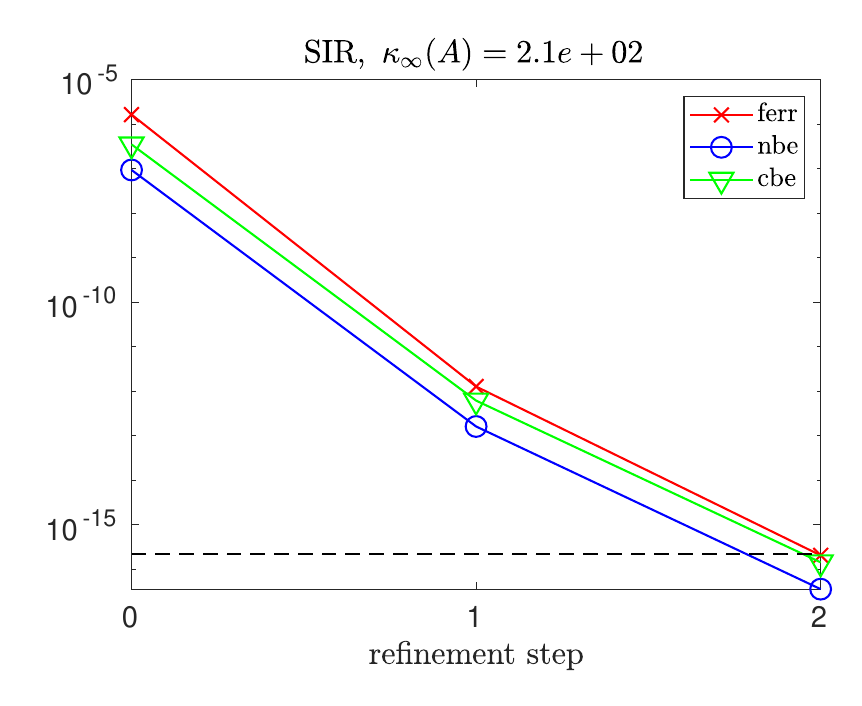}
	\includegraphics[trim={3cm 8cm 4cm 8cm},clip, width=.45\textwidth]{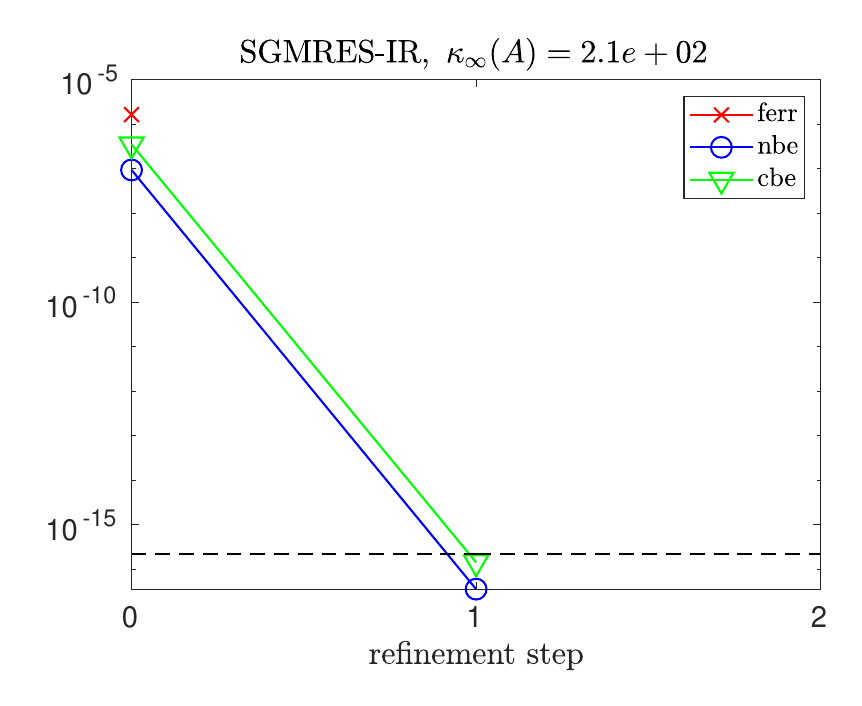}\\
	\includegraphics[trim={3cm 8cm 4cm 8cm},clip, width=.45\textwidth]{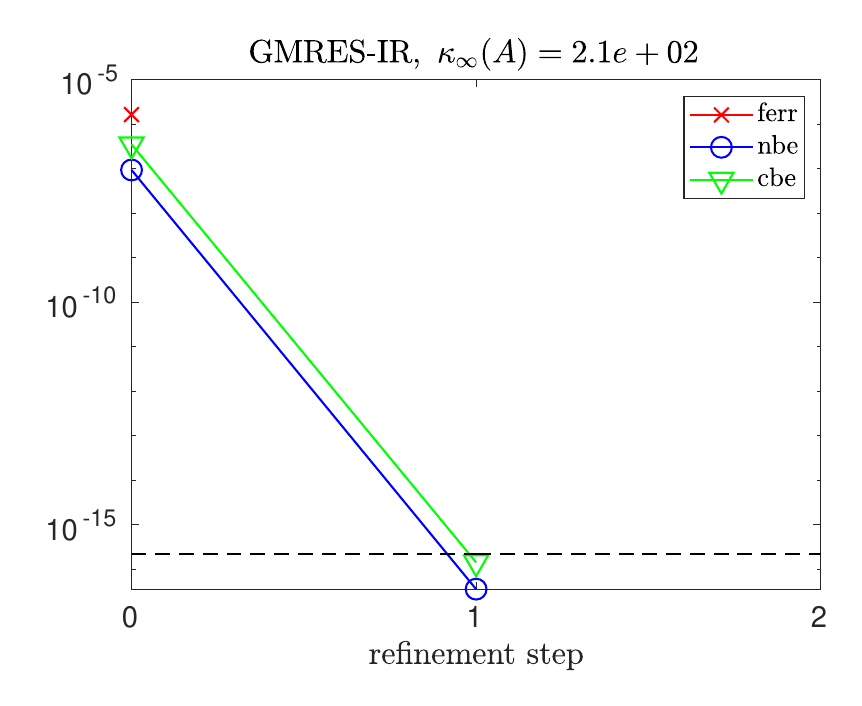}
	\includegraphics[trim={3cm 8cm 4cm 8cm},clip, width=.45\textwidth]{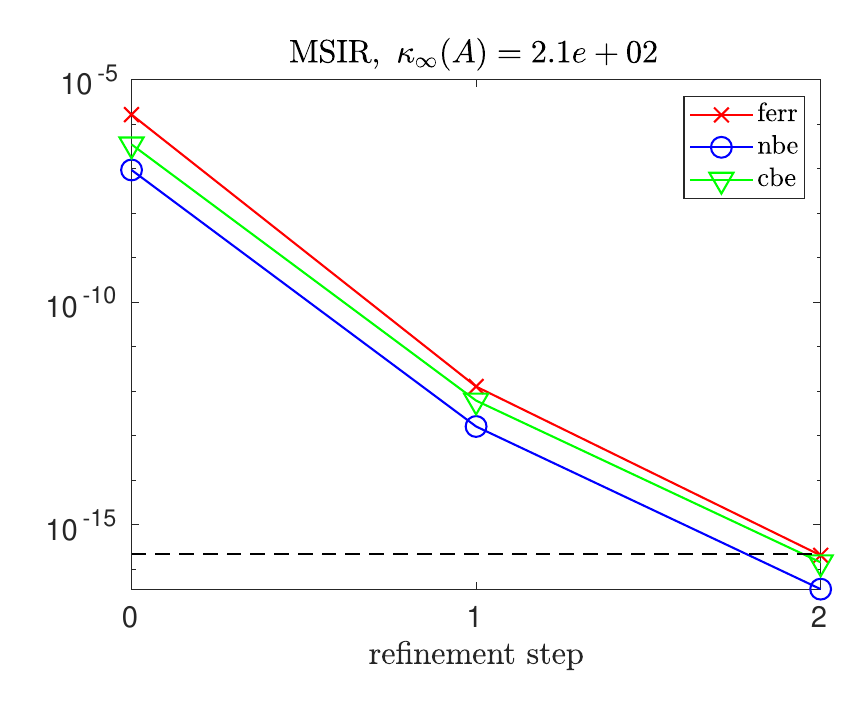}
	\caption{Convergence of errors for a $100 \times 100$ random dense matrix with $\kappa_2 (A)=10^1$ using SIR (top left), SGMRES-IR (top right), GMRES-IR (bottom left), and MSIR (bottom right), with initial precisions $(u_f,u,u_r)$ = (single, double, quad). Note that for SGMRES-IR and GMRES-IR, the forward error is measured as 0 (in double precision) after one refinement step. }
	\label{fig:intro_e1}

\end{figure}

\begin{figure}
\centering
	\includegraphics[trim={3cm 8cm 4cm 8cm},clip, width=.45\textwidth]{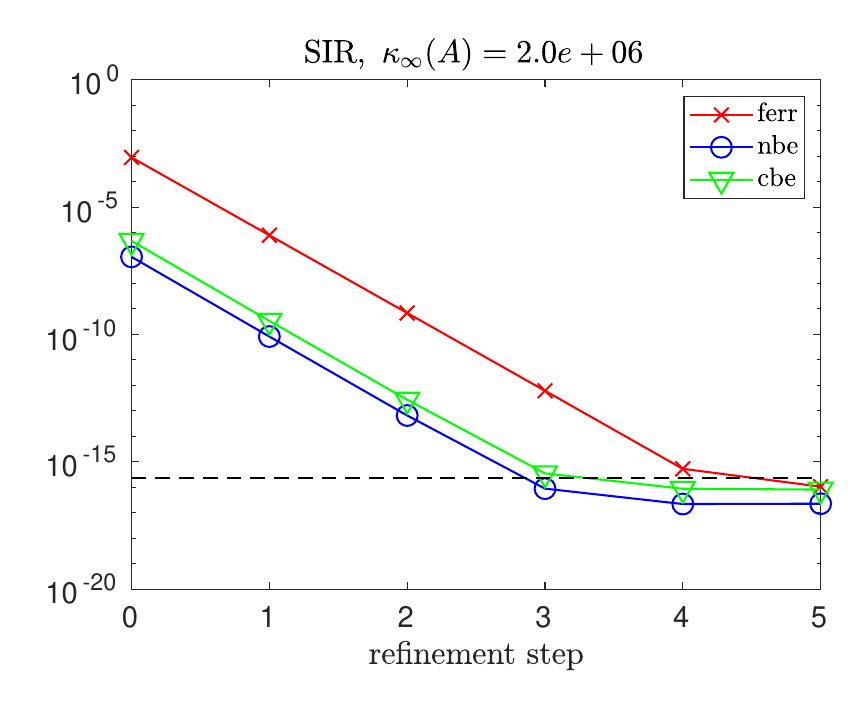}
	\includegraphics[trim={3cm 8cm 4cm 8cm},clip, width=.45\textwidth]{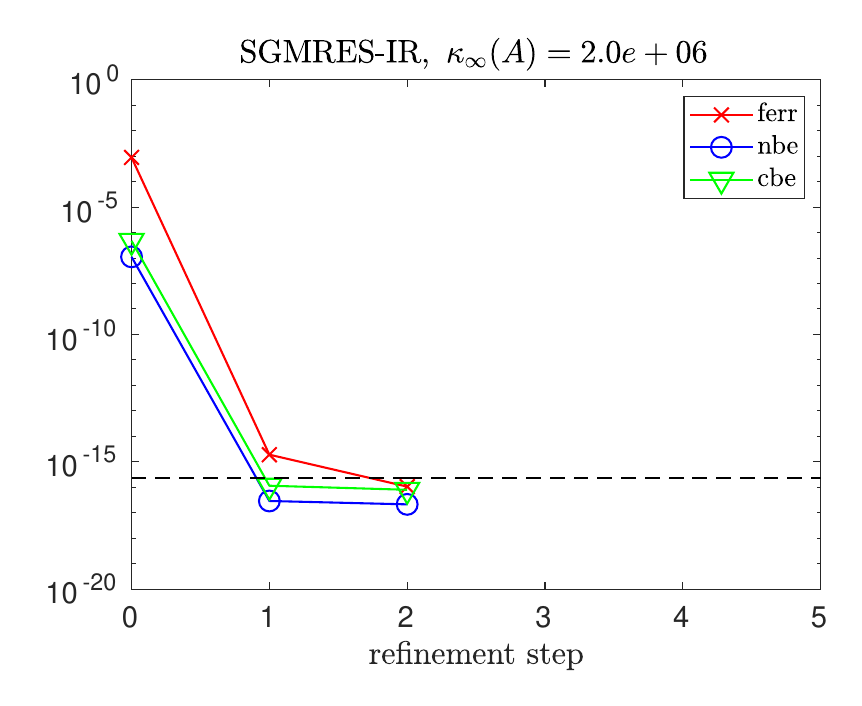}\\
	\includegraphics[trim={3cm 8cm 4cm 8cm},clip, width=.45\textwidth]{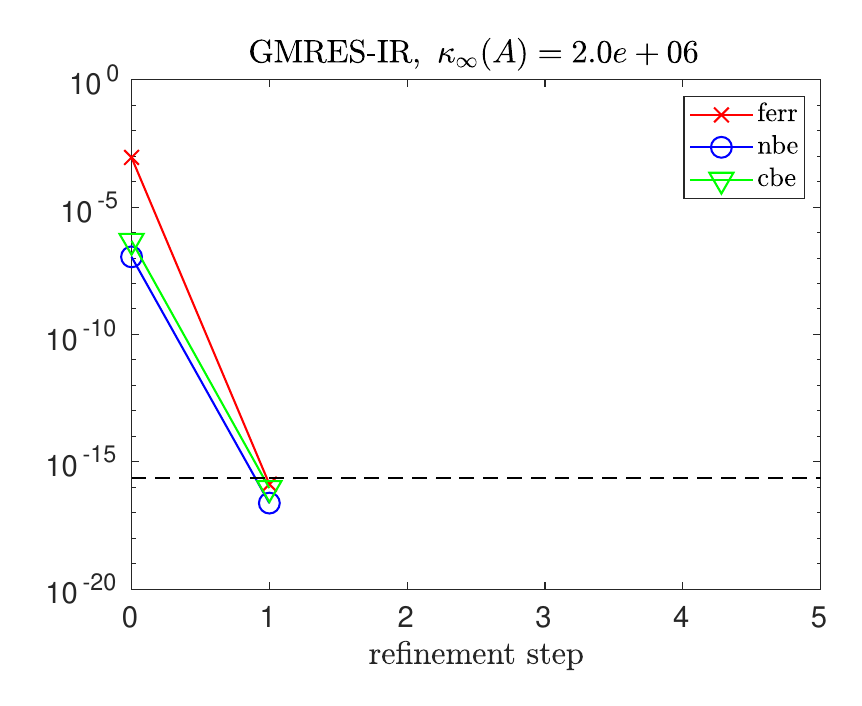}
	\includegraphics[trim={3cm 8cm 4cm 8cm},clip, width=.45\textwidth]{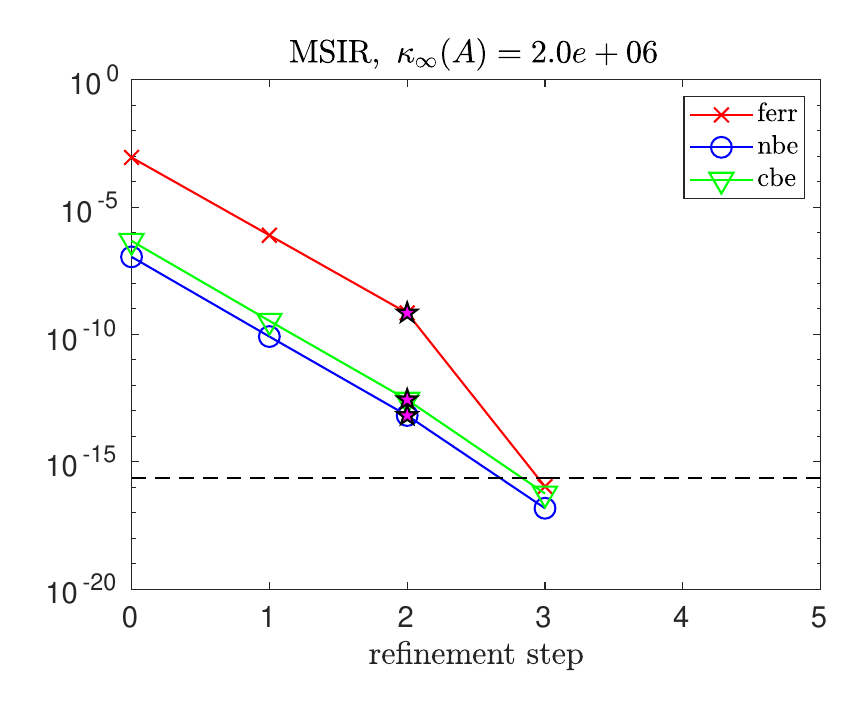}
	\caption{Convergence of errors for a $100 \times 100$ random dense matrix with $\kappa_2 (A)=10^5$ using SIR (top left), SGMRES-IR (top right), GMRES-IR (bottom left), and MSIR (bottom right), with initial precisions $(u_f,u,u_r)$ = (single, double, quad).}
	\label{fig:intro_e5}
\end{figure}

\begin{figure}
\centering
	\includegraphics[trim={3cm 8cm 4cm 8cm},clip, width=.45\textwidth]{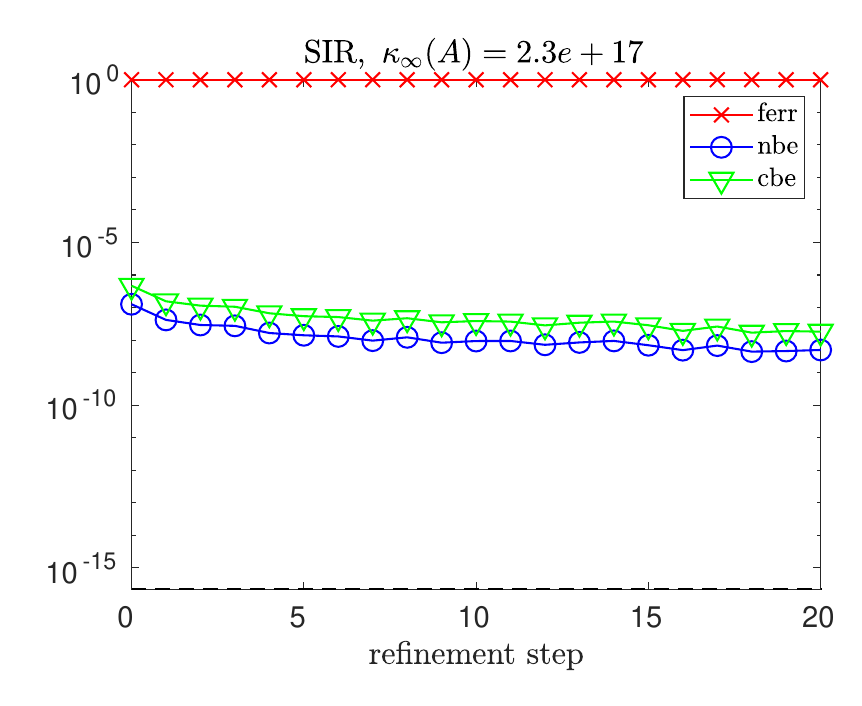}
	\includegraphics[trim={3cm 8cm 4cm 8cm},clip, width=.45\textwidth]{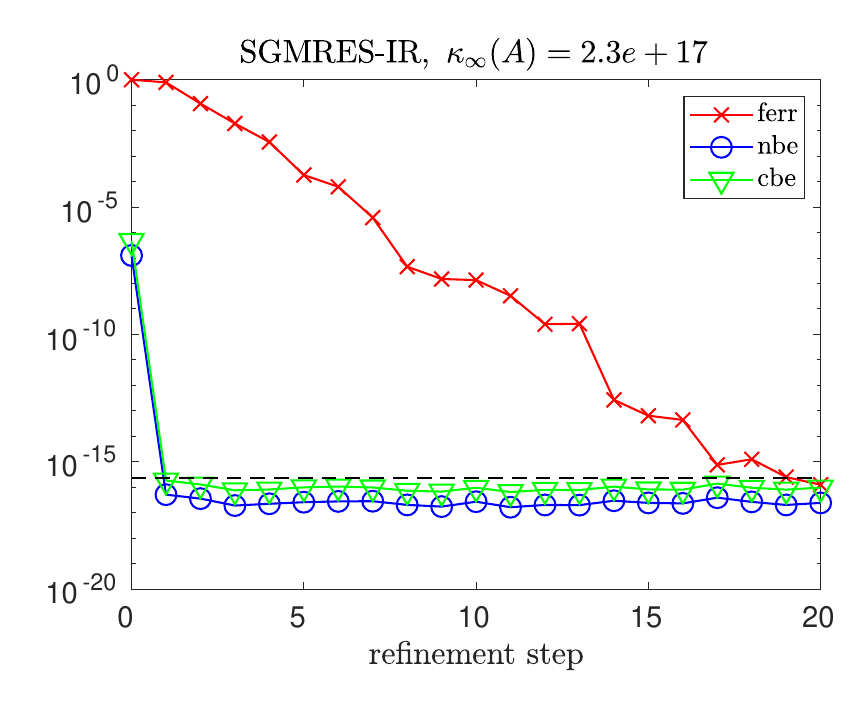}\\
	\includegraphics[trim={3cm 8cm 4cm 8cm},clip, width=.45\textwidth]{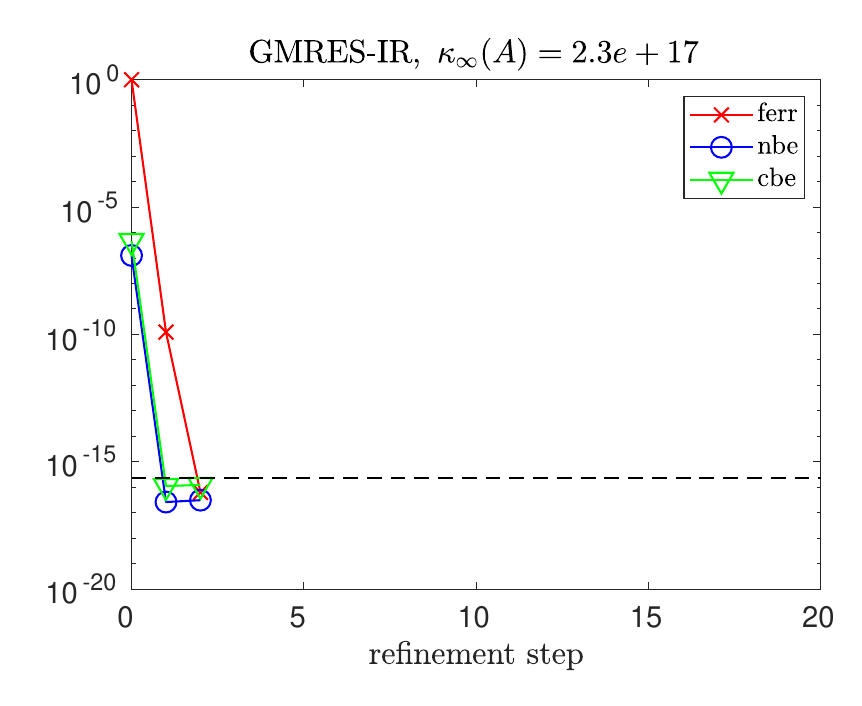}
	\includegraphics[trim={3cm 8cm 4cm 8cm},clip, width=.45\textwidth]{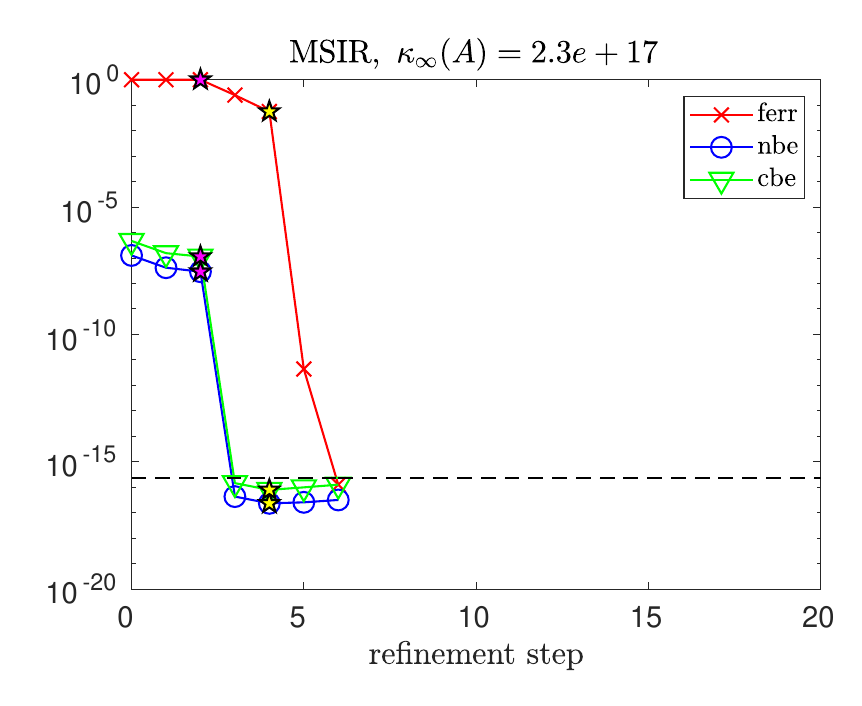}
	\caption{Convergence of errors for a $100 \times 100$ random dense matrix with $\kappa_2 (A)=10^{16}$ using SIR (top left), SGMRES-IR (top right), GMRES-IR (bottom left), and MSIR (bottom right), with initial precisions $(u_f,u,u_r)$ = (single, double, quad).}
	\label{fig:intro_e16}
\end{figure}

Figures \ref{fig:intro_e1}, \ref{fig:intro_e5}, and \ref{fig:intro_e16} show results for condition numbers $10^1, 10^5$, and $10^{16}$, respectively. The red, blue, and green lines in the figures show the behavior of the forward error \verb|ferr| (red), normwise relative backward error \verb|nbe| (blue), and componentwise relative backward error \verb|cbe| (green). The dotted black line shows the value of the initial working precision $u$. If there is a switch in MSIR, the step at which it happened is marked with a star. For instance, if MSIR uses both SGMRES-IR and GMRES-IR, then there are two stars: the first one marks the switch from SIR to SGMRES-IR and the second one marks the switch from SGMRES-IR to GMRES-IR.

For figures in this study, forward error is calculated as $\|\hat{x}-x\|_\infty/\|x\|_\infty$. The normwise relative backward error for $\hat{x}$ is calculated as
\[
\dfrac{\|b-A\hat{x}\|}{\|A\|\|\hat{x}\|+\|b\|},
\]whereas for the computation of the componentwise relative backward error,
\[
\max_i\dfrac{|b-A\hat{x}|_i}{(|A||\hat{x}|+|b|)_i}
\]is used.

The number of refinement steps performed by SIR, SGMRES-IR, GMRES-IR, and MSIR for these matrices are presented in Table \ref{tab:intro}. For SIR, the number in each row gives the number of refinement steps. For SGMRES-IR and GMRES-IR, the number of parenthetical elements gives the total number of refinement steps, and element $i$ in the list gives the number of GMRES iterations performed in refinement step $i$. For example, $(3,4)$ indicates that 2 refinement steps were performed, the first of which took 3 GMRES iterations and the second of which took 4. For the MSIR column, the data for SIR, SGMRES-IR, and GMRES-IR is comma-separated. For example, $2, (2)$ indicates that 2 SIR steps were performed, the algorithm then switched to SGMRES-IR, and performed 1 SGMRES-IR step which required 2 GMRES iterations. Since there is no second set of parentheses, this indicates that there was no switch to GMRES-IR. In all columns, a dash denotes that the algorithm diverged or made no progress; to enable a fair comparison, in our experiments we set $i_{max}$ to be a very high value (here $i_{max}=2000$) in order to allow all approaches that eventually converge sufficient time to do so.


\begin{table}[h!]
	\centering
	\caption{Number of SIR, SGMRES-IR, GMRES-IR, and MSIR steps with the number of GMRES iterations for each SGMRES-IR and GMRES-IR step for initial precisions $(u_f,u,u_r)$ = (single, double, quad).}
	\begin{tabular}{|cc|cccc|}
		\hline
		$\kappa_\infty(A)$ & $\kappa_2(A)$            & \multicolumn{1}{c}{SIR} & \multicolumn{1}{c}{SGMRES-IR} & \multicolumn{1}{c}{GMRES-IR} & MSIR      \\ \hline
		 $2\cdot 10^2$ &$10^1$  & 2                        & (2)                             & (2)                   & 2         \\
		$2 \cdot 10^6$ & $10^5$   & 5                        & (2,3)                              & (2)                & 2, (2) \\
		$2\cdot10^{17}$ & $10^{16}$ &  -                        & \footnotesize{(3,3,3,4,4,4,4,8,100,100,50,12,12,5,4,5,4,4,4,4)}                              & (3,4)                & 2, (3,3), (3,4)   \\ \hline
	\end{tabular}
	\label{tab:intro}
\end{table}

From Figure \ref{fig:intro_e1}, we can see that for well-conditioned matrices, SIR quickly converges to the solution. When the condition number increases closer to $u_f^{-1}=(5.96\cdot 10^8)^{-1} \approx 1.7\cdot 10^7$ as in Figure \ref{fig:intro_e5}, however, SIR convergence begins to slow down, and MSIR switches to SGMRES-IR, which then converges in one step. When the matrix becomes very ill conditioned as in Figure \ref{fig:intro_e16}, then SIR diverges, SGMRES-IR converges very slowly, and thus MSIR makes the second switch to GMRES-IR. 

This illustrates the benefit of the multistage approach. In the case that the problem is well conditioned enough that SIR suffices, MSIR will only use SIR. For the $\kappa_2(A)=10^1$ case, MSIR behaves exactly the same as SIR, and is thus less expensive than SGMRES-IR and GMRES-IR since a single GMRES iteration is more expensive than an SIR step. In the extremely ill-conditioned case, both SIR and SGMRES-IR fail to converge or convergence is too slow. MSIR however does converge quickly, although at roughly double the cost of GMRES-IR. This is the inherent tradeoff. Compared to GMRES-IR, we expect MSIR to be less expensive when the problem is well or reasonably well conditioned and thus only SIR or SGMRES-IR are used (see a comparison of costs in Table \ref{tab:cost}). For cases where the problem is extremely ill conditioned relative to the working precision and GMRES-IR converges, we expect MSIR to be a constant factor more expensive than GMRES-IR. Compared to SIR and SGMRES-IR, the benefit is clear: when SIR and/or SGMRES-IR converges reasonably quickly, MSIR will also converge at roughly the same cost, but MSIR can converge for problems where SIR and SGMRES-IR may not. 


\subsection{Algorithm details}

We now discuss the MSIR algorithm (Algorithm \ref{alg:tsir}) in more detail. The algorithm uses $i$ to denote the global number of refinement steps of any type. The \texttt{iter} variable counts the number of refinement steps of the current refinement variant; e.g., it will first count the number of SIR steps, and is then reset to 0 if the algorithm switches to SGMRES-IR. Unlike the single-stage algorithms, here the input parameter $i_{max}$ specifies the maximum number of refinement steps of each type, meaning, e.g., we will perform up to $i_{max}$ SIR steps before switching to SGMRES-IR. 

\subsubsection{Stopping criteria and convergence detection} \label{convcriteria}

In \cite{dh:06}, Demmel et al. analyze an iterative refinement scheme in which extra precision is used in select computations and provide reliable normwise and componentwise error bounds for the computed solution as well as stopping criteria. They devise a new approach that adaptively changes the precision with which the approximate solution is stored based on monitoring the convergence; if consecutive corrections to the approximate solution are not decreasing at a sufficient rate, the precision of the approximate solution is increased. This has the effect of improving the componentwise accuracy. 

Following this strategy for monitoring the behavior of iterative refinement from \cite{dh:06}, the MSIR algorithm will switch to the next variant if any of the following conditions applies:
\begin{enumerate}
	\item $\frac{\|d_{i+1}\|_\infty}{\|x_{i}\|_\infty}\leq u$ (the correction $d_{i+1}$ changes solution $x_{i}$ too little),
	\item $\frac{\|d_{i+1}\|_\infty}{\|d_{i}\|_\infty}\geq \rho_{thresh}$ (convergence slows down sufficiently),
	\item $\texttt{iter} \geq i_{max}$ (too many iterations of a particular variant have been performed), and, 
	\item $k_{GMRES} \geq k_{max}$ (too many GMRES iterations are performed in one step of SGMRES-IR or GMRES-IR)
\end{enumerate}
where $d_{i}$ is the correction of the solution $x_{i}$, $\rho_{thresh}\in \mathbb{R}^+$ is a threshold for convergence and $\rho_{thresh}<1$, $i_{max} \in \mathbb{N}^+$ is the maximum number of iterations, $k_{max} \in \mathbb{N}^+$ 
is the maximum number of GMRES iterations performed per SGMRES-IR step, 
and $u$ is the working precision. 

As is explained in \cite{dh:06}, the analyses of Bowdler in 1966 \cite{b:66} and Moler in 1967 \cite{m:67} showed that $\|d_{i+1}\|_\infty$ should decrease by a factor of at most $\rho = O(u_f)\kappa_\infty(A)$ at each step, and thus the solution $x_{i+1}$ should converge like a geometric sum, meaning that $\sum_{j=i+1}^{\infty} \Vert d_i \Vert_\infty \leq \Vert d_{i+1} \Vert_\infty/(1-\rho)$. This geometric convergence will end when rounding errors become significant. So if we have $\Vert d_{i+1} \Vert_\infty/ \Vert d_i \Vert_\infty \geq \rho_{thresh}$, it either means that (1) convergence has slowed down due to rounding errors, or (2) convergence is slow from the beginning as the quantity $\rho$ is close to 1 (meaning that the problem is too ill conditioned with respect to the precision $u_f$). In the case that $\|d_{i+1}\|_\infty > \|d_{i}\|_\infty$, this indicates that the iterative refinement process is diverging, again because the problem is too ill conditioned with respect to the precision $u_f$.

We reiterate that the fourth condition above only applies to the switch from SGMRES-IR to GMRES-IR, or GMRES-IR to SIR with increased precision(s). We note that as in line \ref{tsir:nancheck} in Algorithm \ref{alg:tsir}, in case SIR produces a correction $d_{i+1}$ containing Inf or NaN, we immediately switch to SGMRES-IR without performing the solution update. This situation can arise as a result of performing the triangular solves in low precision. We have implemented a simple scaling to help avoid this situation, although in certain scenarios it may still occur; see further details in Section \ref{sec:scaling}. We note that this is less of a concern for GMRES-based approaches, since within GMRES the triangular solves are performed in either precision $u$ (SGMRES-IR) or $u^2$ (GMRES-IR).

Moreover, the convergence detection in lines \ref{conv1}, \ref{conv2}, and \ref{conv3} in Algorithm \ref{alg:tsir} can be performed using the normwise relative error estimate discussed in \cite{dh:06},
\begin{equation}
	max\left\lbrace \dfrac{\frac{\|d_{i+1}\|_\infty}{\|x_{i}\|_\infty}}{1-\rho_{max}},\gamma u\right\rbrace \approx \dfrac{\|x_{i}-x\|_\infty}{\|x\|_\infty},
	\label{eq:convdetection}
\end{equation}
where $\gamma = max(10,\sqrt{n})$, $n$ is the size of the matrix $A$, and $\rho_{max} := max_{j \leq i}\frac{\|d_{j+1}\|_\infty}{\|d_{j}\|_\infty}$ is the maximum ratio of successive corrections. The quantity $(\|d_{i+1}\|_\infty/\|x_{i}\|_\infty)/(1-\rho_{max})$ is stored in Algorithm \ref{alg:tsir} as the quantity $\phi_i$. If the condition in lines \ref{conv1}, \ref{conv2}, and \ref{conv3} is triggered, one could then test for convergence by checking whether both $\phi_i \leq \sqrt{n}u$ and $\phi_i \geq 0$, and if so, declare convergence (note that $\phi_i <0$ indicates divergence). If convergence is not detected, we must also decide whether to keep the current approximate solution or reset the approximate solution to $x_0$, since this will be input to the next stage. We test if $\phi_i > \phi_0$, and if so, we reset the initial solution. 
We note that other measures, such as for the componentwise error, could also be used to detect convergence. 

We note that the criteria used to determine whether the current algorithm should quit iterating are one iteration behind; in other words, we can only detect convergence (or non-convergence) in step $i$ after computing the correction term in step $i+1$. For this reason, the MSIR algorithm will generally compute at least two steps of a particular variant before deciding to switch. The two cases where fewer than 2 steps of a variant will occur are 1) when the computed update contains Infs or NaNs (as in line \ref{tsir:nancheck} in Algorithm \ref{alg:tsir}), in which case we switch without updating the current solution, and 2) when, for SGMRES-IR and GMRES-IR, the number of GMRES iterations exceeds $k_{max}$, since this is detectable immediately in the current step.

In \cite{dh:06}, to determine $\rho_{thresh}$, the authors define `cautious' and `aggressive' settings. Cautious settings produce maximally reliable error bounds for well-conditioned problems, whereas aggressive ones will lead to more steps on the hardest problems and usually, but not always, give error bounds within a factor of 100 of the true error. For cautious mode, the authors suggest $\rho_{thresh}$ should be set to $0.5$, which was used by Wilkinson \cite{w:63}, and for aggressive mode, $0.9$. In our experiments we always use the cautious setting, but we note that $\rho_{thresh}$, $i_{max}$, and $k_{max}$ should be set according to the relative costs of the different refinement schemes in practice.

\subsubsection{Scaling}
\label{sec:scaling}
Because of the smaller range of half precision, simply rounding higher precision quantities to a lower precision can cause overflow, underflow, or the introduction of subnormal numbers. 
In \cite{higham2019squeezing}, the authors develop a new algorithm for converting single and double precision quantities to half precision. This algorithm involves performing a two-sided scaling for equilibration and then an additional scaling is performed to make full use of the range of half precision before finally rounding to half precision. In our algorithm, when half precision is used for the factorization, we first attempt an LU factorization without scaling. We then test whether the resulting L and U factors contain Inf or NaN; if so (marked with a * in the tables in Section \ref{sec:results}), we retry the LU factorization in line \ref{tsir:lufact} using the two-sided scaling algorithm of \cite{higham2019squeezing}. In all cases, regardless of what precision is used for the factorization, after computing $x_0$ in line \ref{tsir:initsolve} we test whether the initial solution contains Inf or NaN; if so, we simply use the zero vector as the initial approximate solution and proceed. 

We also incorporate scaling in each refinement step. After computing the residual $r_i$, we scale the result to obtain $r_i = r_i/\Vert r_i\Vert_{\infty}$ (line \ref{r1} in Algorithm \ref{alg:tsir}). This scaling is then undone when we update the approximate solution, via $x_{i+1}=x_i + \Vert r_i\Vert_{\infty} d_{i+1}$ (lines \ref{x1}, \ref{x2}, and \ref{x3} in Algorithm \ref{alg:tsir}). As long as $1/\Vert A\Vert_{\infty}$ does not underflow and $\Vert A^{-1} \Vert_{\infty}$ does not overflow, then this scaling avoids the largest element of $d_{i+1}$ overflowing or underflowing. 


\subsection{Error bounds for different variants}
\label{sec:errorbounds}

There are various scenarios on the usage of precisions which will yield different error bounds and different constraints on condition number. Besides three-precision variants, two precisions or a fixed precision can be used in the algorithms discussed in this study. We summarize the convergence criteria for the precision combinations used in our approach and refer the reader to \cite{ch:18} and \cite{h:21} for more general bounds. In particular, we assume that $u_f \geq u$ and $u_r \leq u^2$. We also restrict ourselves to IEEE precisions (see Table \ref{tab:eps}), although we note that alternative formats like bfloat16 \cite{bfloat16} could also be used. 

Under the assumptions that $u_f \geq u$ and $u_r \leq u^2$, for SIR, both the relative forward and backward errors are guaranteed to converge to the level of the working precision when $\kappa_\infty(A)<u_f^{-1}$. For SGMRES-IR and GMRES-IR, the constraints for convergence of forward and backward errors differ. Summarizing the analysis in \cite{h:21}, for our particular SGMRES-IR variant, the constraint on convergence of the backward error is $\kappa_\infty(A)<u^{-1/3}u_f^{-1/3}$ and the constraint on the convergence of the forward error is $\kappa_\infty(A)<u^{-1/3}u_f^{-2/3}$. It is interesting to note that, as an artifact of the analysis, the constraint for convergence of the backward error for this variant is more strict than that for the forward error. However, since the backward error is bounded by the forward error, the constraint for both backward and forward error to converge to the level $u$ can be given as $\kappa_\infty(A)<u^{-1/3}u_f^{-2/3}$. 

In \cite{ch:18}, the constraint for convergence of the backward error to level $u$ in GMRES-IR given as is $\kappa_\infty(A)<u^{-1}$. However, the authors in \cite{h:21} point out that this bound relies on assumptions that may be overly optimistic, and revise this to the tighter constraint $\kappa_\infty(A)<u^{-1/2}u_f^{-1/2}$. The constraint for the convergence of the forward error to this level in GMRES-IR given in \cite{ch:18} is $\kappa_\infty(A)<u^{-1/2}u_f^{-1}$. Thus the tighter constraint for the convergence of backward error in \cite{h:21} is again more strict than the constraint for the convergence of the forward error, and since the backward error is bounded by the forward error, we can take $\kappa_\infty(A)<u^{-1/2}u_f^{-1}$ to be the constraint for the convergence of both backward and forward error to level $u$ in GMRES-IR.

In Table \ref{tab:kappalimit} we quantify these constraints on $\kappa_\infty(A)$ required for convergence of the normwise relative backward and forward errors to the level of the working precision for various precision combinations. To compute the constraints on $\kappa_\infty(A)$, we have used the unit roundoff values in Table \ref{tab:eps}.

\begin{table}[]
\centering
\caption{Various IEEE precisions and their units roundoff. }
\label{tab:eps}
\begin{tabular}{|l|l|}
\hline
\multicolumn{1}{|c|}{Precision} & \multicolumn{1}{c|}{Unit Roundoff} \\ \hline
fp16 (half) & $4.88\cdot 10^{-4}$ \\ \hline
fp32 (single) & $5.96\cdot 10^{-8}$ \\ \hline
fp64 (double) & $1.11\cdot 10^{-16}$ \\ \hline
fp128 (quad) & $9.63\cdot 10^{-35}$ \\ \hline
\end{tabular}
\end{table}

\begin{table}[]
\centering
\caption{Constraints on $\kappa_\infty(A)$ for which the relative forward and normwise backward errors are guarantee to converge to the level $u$ for a given combination of precisions for the different variants of IR. }
\label{tab:kappalimit}
\begin{tabular}{ccc|ccc}
$u_f$ & $u$ & $u_r$ & SIR & SGMRES-IR & GMRES-IR \\ \hline
half & single & double & $2\cdot 10^{3}$ & $4\cdot 10^{4}$ & $8\cdot 10^{6}$ \\
single & single & double & $2\cdot 10^{7}$ & $2\cdot 10^{7}$ & $7\cdot 10^{10}$\\
half & double & quad & $2\cdot 10^{3}$ & $3\cdot 10^{7}$ & $2\cdot 10^{11}$ \\
single & double & quad & $2\cdot 10^{7}$ & $1\cdot 10^{10}$ & $2\cdot 10^{15}$ \\
double & double & quad & $9\cdot 10^{15}$ & $9\cdot 10^{15}$ & $9\cdot 10^{23}$
\end{tabular}
\end{table}

\section{Numerical experiments}\label{sec:results}

In this section we present numerical experiments performed in MATLAB on a number of synthetic and real-world matrices from SuiteSparse \cite{dh:11} for MSIR and single-stage iterative refinement variants in three precisions. The problems we test are small and are meant to demonstrate the numerical behavior of the algorithms. Performance results for SIR and GMRES-based variants for larger problems on modern GPUs can be found in, e.g., \cite{ht:18}. 

For quadruple precision, we use the Advanpix multiprecision computing toolbox for MATLAB \cite{advanpix}. To simulate half precision floating point arithmetic, we use the \texttt{chop} library and associated functions from \cite{higham2019simulating}, available in the repositories \texttt{https://github.com/higham/chop} and \texttt{https://github.com/SrikaraPranesh/LowPrecision\_Simulation}. For single and double precision we use the built-in MATLAB datatypes. All experiments were performed on a computer with an Intel Core i7-9750H processor with 12 CPUs and 16 GB RAM with OS system Ubuntu 20.04.1 using MATLAB 2020a. Our MSIR algorithm and associated functions are available through the repository \texttt{https://github.com/edoktay/msir}, which includes scripts for generating the data and plots in this work.

In all experiments, we use $i_{max}=2000$ and $\rho_{thresh}=0.5$. We have chosen to use a very high value for $i_{max}$ in both MSIR and single-stage algorithms as it allows us to see the true convergence behavior of each algorithm (and the benefits of MSIR in cases where convergence for single-stage algorithms happens eventually but is very slow). Of course, in practice, one would use a much smaller value.

For the GMRES convergence tolerance, we use $\tau = 10^{-6}$ when the working precision is single and $\tau = 10^{-10}$ when the working precision is double. We set $k_{max} = 0.1n$. For purposes of discerning the actual attainable accuracy of MSIR, in the experiments we explicitly compute and plot the computed forward and backward errors and use these as stopping criteria rather than the error estimate given by \eqref{eq:convdetection}. For a fair numerical comparison, we also apply the scalings discussed in Section \ref{sec:scaling} in SIR, SGMRES-IR, and GMRES-IR. 

For each variant, we compare the number of refinement steps required for forward and backward errors to reach the level of accuracy corresponding to the initial working precision. For GMRES-based approaches and MSIR, the number of GMRES iterations per step is given parenthetically (see the explanation in Section \ref{sec:tsir}). For the MSIR results, a semicolon indicates a precision switch: the factorization precision is doubled (and other precisions increased if necessary to ensure $u_f \geq u$ and $u_r \leq u^2$) and the algorithm restarts with SIR. A dash (-) in the tables indicates that forward and/or backward errors are not decreasing. Select convergence plots for MSIR are presented in figures while the number of steps and iterations for all variants are shown in tables.


Again, the red, blue, and green lines in the figures show the behavior of the forward
error \verb|ferr| (red), normwise relative backward error \verb|nbe| (blue), and componentwise relative backward error \verb|cbe| (green). The dotted black line shows the value of the initial working precision $u$. Switches in MSIR are denoted by stars. A magenta star indicates a switch from SIR to SGMRES-IR, a yellow start indicates a switch from SGMRES-IR to GMRES-IR, and a cyan star indicates switch from GMRES-IR to SIR with increased precision(s).

\subsection{Random dense matrices} \label{sec:randnmat}

We first test our algorithm on random dense matrices of dimension $n=100$ generated in MATLAB via the command \verb|gallery('randsvd',n,kappa(i),2)| (Section \ref{sec:randnmat_mode2}) and the command \verb|gallery('randsvd',n,kappa(i),3)| (Section \ref{sec:randnmat_mode3}), where \verb|kappa| is the array of the desired 2-norm condition numbers $\kappa_2 (A) =$\{$10^1$, $10^2$, $10^4$, $10^5$, $10^7$, $10^9$, $10^{11}$, $10^{14}$\}, and \verb|2| and \verb|3| stand for the modes which dictate the singular value distribution. Mode 2 generates matrices with one singular value equal to $1/\kappa_2(A)$ and the rest of the singular values equal to 1. Mode 3 generates matrices with geometrically distributed singular values. Mode 3 represents a difficult case for GMRES-based iterative refinement methods. It is known that a cluster of singular values near zero causes stagnation of GMRES, and a low precision preconditioner may fail to effectively shift this cluster away from the origin. We test the algorithms using initial precision combinations $(u_f,u,u_r)$ = (single, double, quad), $(u_f,u,u_r)$ = (half, single, double), and $(u_f,u,u_r)$ = (half, double, quad). For simplicity, $b$ is chosen to be a vector of normally distributed numbers generated by the MATLAB command \verb|randn| for all experiments in this section. For reproducibility, we use the MATLAB command \verb|rng(1)| before generating each linear system.

\subsubsection{Random dense matrices with one small singular value}\label{sec:randnmat_mode2}
\begin{figure}[h!]
	\centering
	\includegraphics[trim={3cm 8cm 4cm 8cm},clip, width=.45\textwidth]{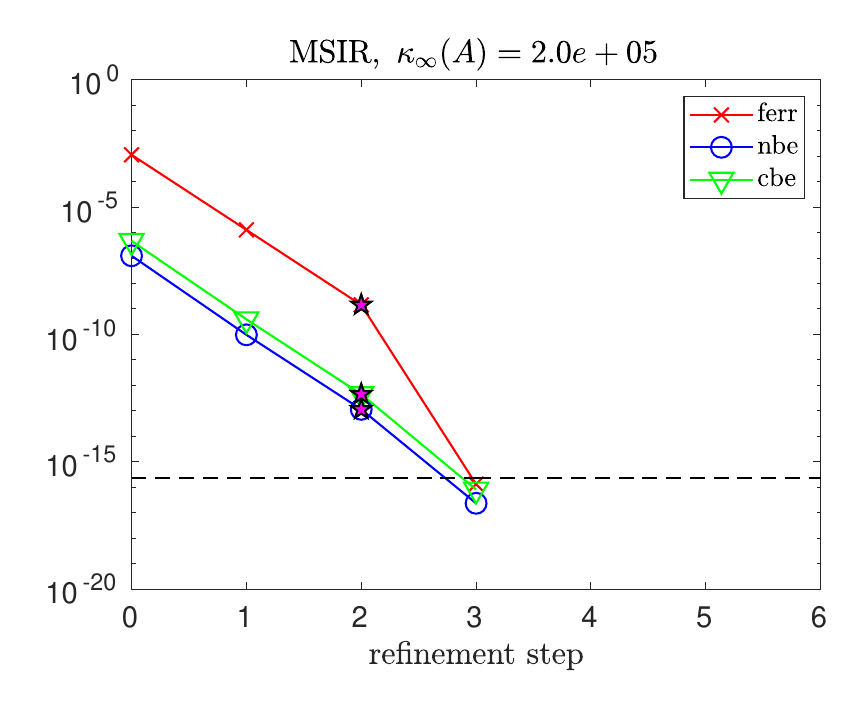}
	\includegraphics[trim={3cm 8cm 4cm 8cm},clip, width=.45\textwidth]{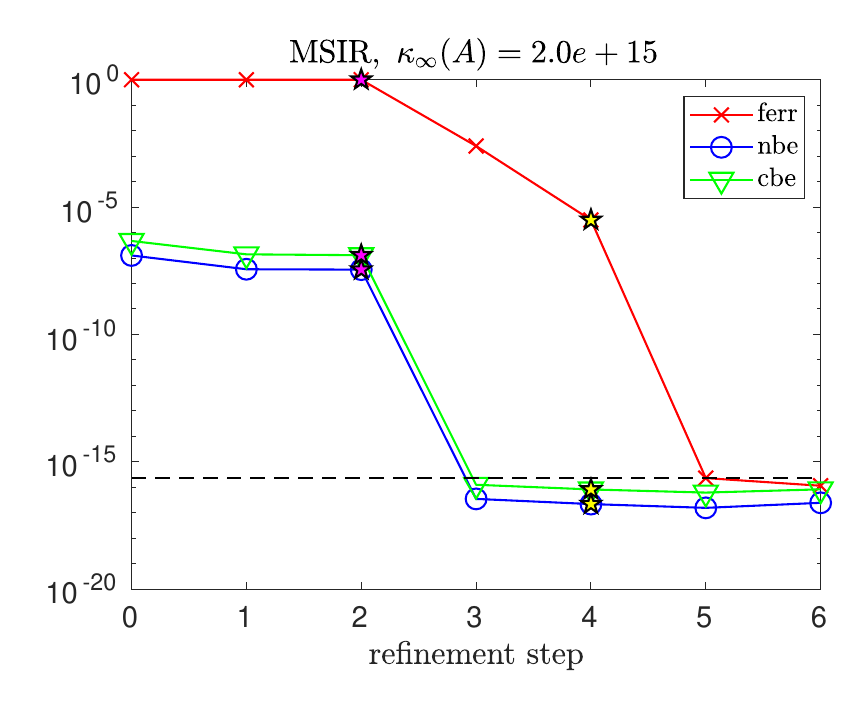}
	\caption{Convergence of errors in MSIR for random dense matrices (mode 2) with $\kappa_2 (A)=10^4$ (left), and $\kappa_2 (A)=10^{14}$ (right) for initial precisions $(u_f,u,u_r)$ = (single, double, quad); see also Table \ref{tab:rand_sdq}.}
	\label{fig:randnmat_124}
\end{figure}

Table \ref{tab:rand_sdq} shows the convergence behavior of SIR, SGMRES-IR, GMRES-IR, and MSIR for matrices generated via the MATLAB command \verb|gallery('randsvd',n,kappa(i),2)| using initial precisions $(u_f,u,u_r)$ = (single, double, quad). 
From Table \ref{tab:rand_sdq}, we can see that SIR no longer converges once $\kappa_\infty(A)$ exceeds $10^{11}$. As $\kappa(A)$ increases, the convergence rate of SIR slows down. At the extreme case of $\kappa_2(A)=10^9$, SIR requires 205 steps to converge! In this case, MSIR thus demonstrates a significant improvement over SIR.
Even in cases where only 5 SIR steps are required, notice that the MSIR algorithm deems this too slow and begins switching to SGMRES-IR after 2 SIR steps. SGMRES-IR converges in just a few refinement steps with a small numbers of GMRES iterations per step except for the most ill-conditioned problem. For the most ill-conditioned matrix, SGMRES-IR alone requires 7 steps to converge. MSIR in this case switches twice: from SIR to SGMRES-IR, and then finally to GMRES-IR. The same is observed when $\kappa_2(A)=10^{11}$. We plot the convergence trajectories of MSIR for the $\kappa_2(A)=10^4$ and $\kappa_2(A)=10^{14}$ problems in Figure \ref{fig:randnmat_124}.

\begin{table}[h!]
	\centering
	\caption{Number of SIR, SGMRES-IR, GMRES-IR, and MSIR steps with the number of GMRES iterations for each SGMRES-IR and GMRES-IR step for random dense matrices (mode 2) with various condition numbers $\kappa_\infty(A)$, $\kappa_2(A)$, using initial precisions $(u_f,u,u_r)$ = (single, double, quad).}
	\begin{tabular}{|cc|cccc|}
		\hline
		$\kappa_\infty (A)$     & $\kappa_2(A)$  & {SIR} & {SGMRES-IR} & {GMRES-IR} & {MSIR} \\  \hline
		$2\cdot 10^2$           & $10^1$           & 2                    & (2)                     & (2)                     & 2                     \\
		$2\cdot 10^3$           & $10^2$           & 3                    & (2)                      & (2)                     & 2, (2)                     \\
		$2\cdot 10^5$           & $10^4$           & 5                    & (2,3)                     & (2)                     & 2, (2)                 \\
		$2\cdot 10^6$           & $10^5$           & 5                    & (2,3)                      & (2)                     & 2, (2)                \\
		$2\cdot 10^8$           & $10^7$           & 19                  & (2,3)                    & (2,3)                  & 2, (2,3)               \\
		$2\cdot 10^{10}$        & $10^9$           & 205                    & (2,2)                      & (2,3)                     & 2, (2,3)               \\
		$2\cdot 10^{12}$        & $10^{11}$        & -                    & (3,3,4)                   & (3,4)                     & 2, (3,3), (3)             \\
		$2\cdot 10^{15}$        & $10^{14}$        & -                    & (3,3,3,4,4,4,4)              & (3,4)       & 2, (3,3), (3,4)        \\ \hline
	\end{tabular}
	\label{tab:rand_sdq}
\end{table}

In Table \ref{tab:rand_hsd}, we show data for experiments with $(u_f,u,u_r)$ = (half, single, double). In this case, as predicted in Table \ref{tab:kappalimit}, SIR fails to converge once $\kappa_\infty(A)>2\cdot 10^5$. Moreover, it is seen that MSIR switches to SGMRES-IR even for relatively well-conditioned matrices for which SIR still converges (i.e., for $\kappa_2(A)=10^2$). Although this switch is technically not necessary, we argue that the difference in cost versus SIR will be a constant factor (SIR requires 3 LU solves in precision $u_f$ whereas here SGMRES-IR requires 2 LU solves in precision $u_f$ and 3 in precision $u$; see Table \ref{tab:cost}). 
SGMRES-IR works well up to condition number $\kappa_\infty(A)=2\cdot 10^6$, and thus MSIR only performs one switch in these cases. 

It is interesting to notice that the theory from \cite{h:21} says that we can only expect SGMRES-IR to converge as long as $\kappa_\infty(A)<4\cdot 10^4$; see also Table \ref{tab:kappalimit}. However, as this and other experiments in this work show, SGMRES-IR as well as GMRES-IR actually often perform beyond the limits given by the analysis. This confirms our motivation for the multistage MSIR approach; we cannot rely entirely on the existing theoretical bounds to determine whether or not an algorithm will be effective. 

\begin{figure}[h!]
	\centering
	\includegraphics[trim={3cm 8cm 4cm 8cm},clip, width=.45\textwidth]{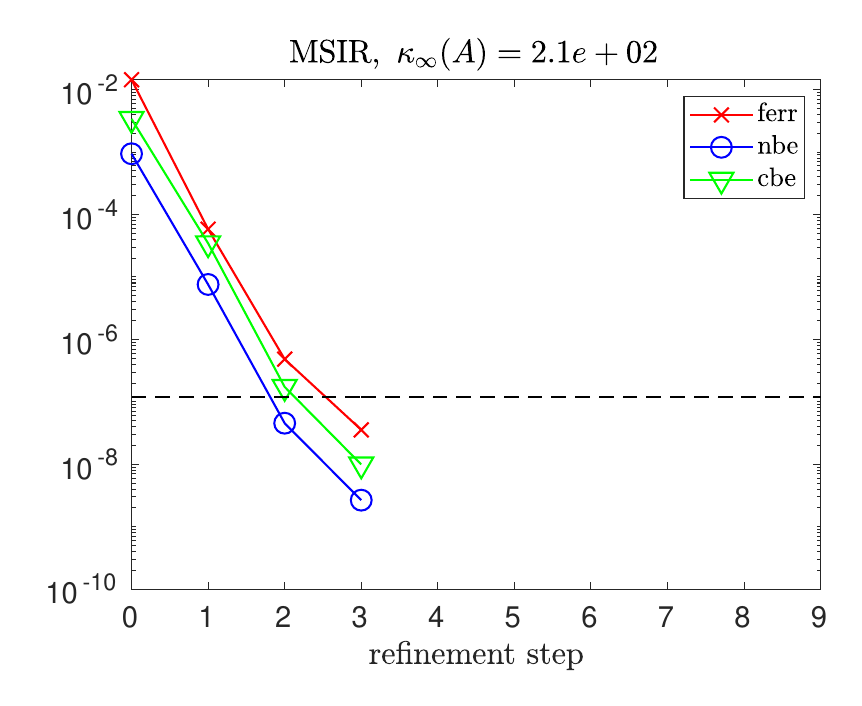}
	\includegraphics[trim={3cm 8cm 4cm 8cm},clip, width=.45\textwidth]{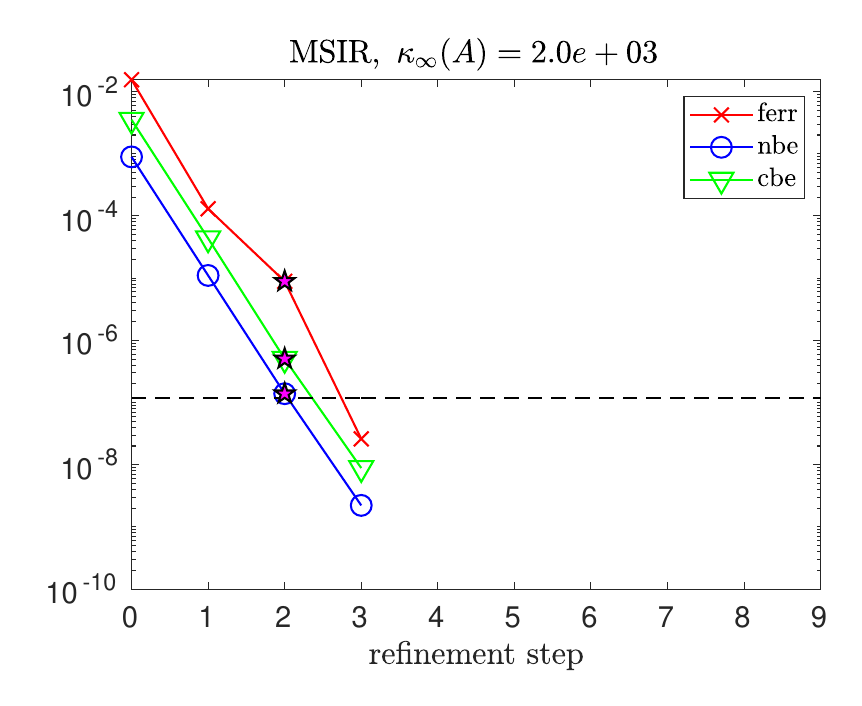}\\
	\includegraphics[trim={3cm 8cm 4cm 8cm},clip, width=.45\textwidth]{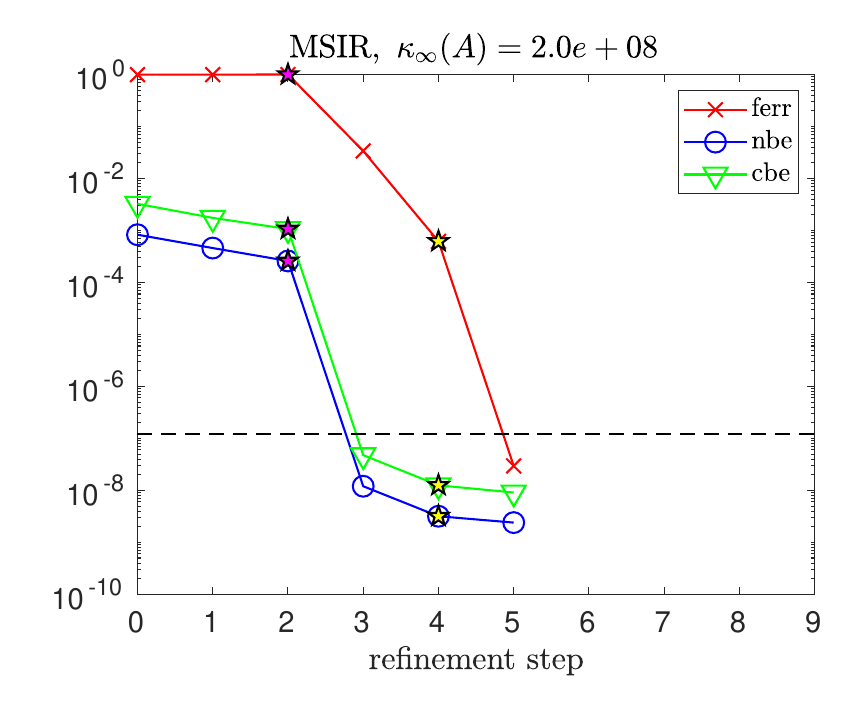}
	\includegraphics[trim={3cm 8cm 4cm 8cm},clip, width=.45\textwidth]{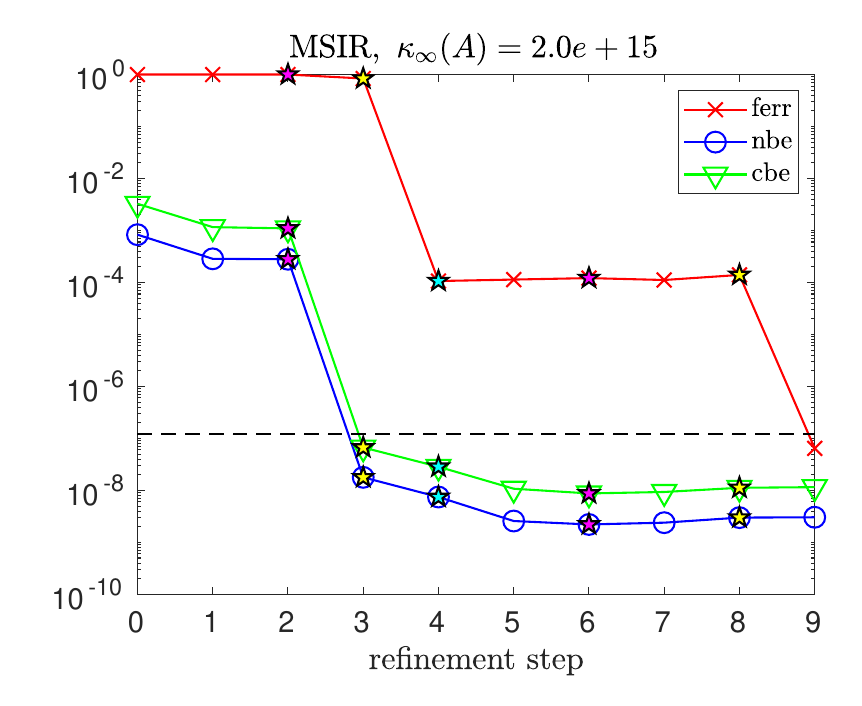}
	\caption{Convergence of errors in MSIR for random dense matrices (mode 2) with $\kappa_2 (A)=10^1$ (top left), $\kappa_2 (A)=10^2$ (top right), $\kappa_2 (A)=10^7$ (bottom left), and $\kappa_2 (A)=10^{14}$ (bottom right) for initial precisions $(u_f,u,u_r)$ = (half, single, double); see also Table \ref{tab:rand_hsd}.}
	\label{fig:randnmat_012}
\end{figure}

For the condition number $\kappa_\infty(A) = 2\cdot 10^8$, SGMRES-IR convergence eventually slows down, and as a result MSIR switches a second time to GMRES-IR. Beyond this limit, for $\kappa_\infty(A)=2\cdot 10^8$, GMRES-IR still converges (despite that the theory only guarantees that GMRES-IR will work up to condition number $8\cdot 10^6$). Starting at $\kappa_\infty(A)=2\cdot 10^{10}$, GMRES convergence slows down significantly (requiring $n$ iterations per refinement step), and thus MSIR switches to precisions 
$(u_f, u, u_r)$ = (single, single, double) and starts again with SIR. We show plots of MSIR convergence behavior for the cases $\kappa_2(A)\in \{10^1, 10^2, 10^7, 10^{14}\}$ in Figure \ref{fig:randnmat_012}. 

\begin{table}[h!]
	\centering
	\caption{Number of SIR, SGMRES-IR, GMRES-IR, and MSIR steps with the number of GMRES iterations for each SGMRES-IR and GMRES-IR step for random dense matrices (mode 2) with various condition numbers $\kappa_\infty(A)$, $\kappa_2(A)$, using initial precisions $(u_f,u,u_r)$ = (half, single, double).}
	\scalebox{0.80}{\begin{tabular}{|cc|cccc|}
		\hline
		 $\kappa_\infty (A)$     & $\kappa_2(A)$ & {SIR} & {SGMRES-IR} & {GMRES-IR} & {MSIR} \\ \hline
		$2\cdot 10^2$           & $10^1$           & 3                    & (3)                        & (3)                       & 3                     \\
		$2\cdot 10^3$           & $10^2$           & 3                    & (3)                        & (3)                       & 2, (3)                 \\
		$2\cdot 10^5$           & $10^4$           & 19                    & (3,4)                      & (3,4)                       & 2, (3,4)               \\
		$2\cdot 10^6$           & $10^5$           & -                    & (3,4)                    & (3,4)                     & 2, (3,4)               \\
		$2\cdot 10^8$           & $10^7$           & -                    & (5,12,9,17,5,5,5)         & (5,5)                     & 2, (5,10), (10)          \\
		$2\cdot 10^{10}$        & $10^9$           & -                    & -                          & (100,100,100,22,11,3)              & 2, (5,9), (10); 2, (2,2,2), (2,3,3)                     \\
		$2\cdot 10^{12}$        & $10^{11}$        & -                    & -                          & (100,63,10)                 & 2, (10), (10); 2, (2,2), (2)                     \\
		$2\cdot 10^{15}$        & $10^{14}$        & -                    & -                          & (100,63,10)                & 2, (10), (10); 2, (2,2), (2)                     \\ \hline
	\end{tabular}}
	\label{tab:rand_hsd}
\end{table}

\begin{figure}[h!]
	\centering
	\includegraphics[trim={3cm 8cm 4cm 8cm},clip, width=.45\textwidth]{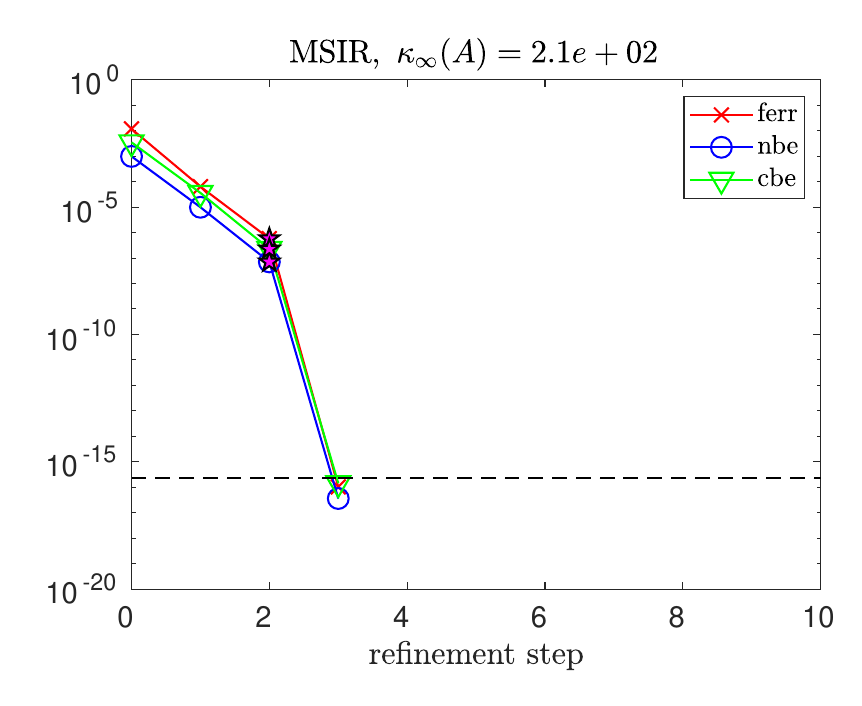}
	\includegraphics[trim={3cm 8cm 4cm 8cm},clip, width=.45\textwidth]{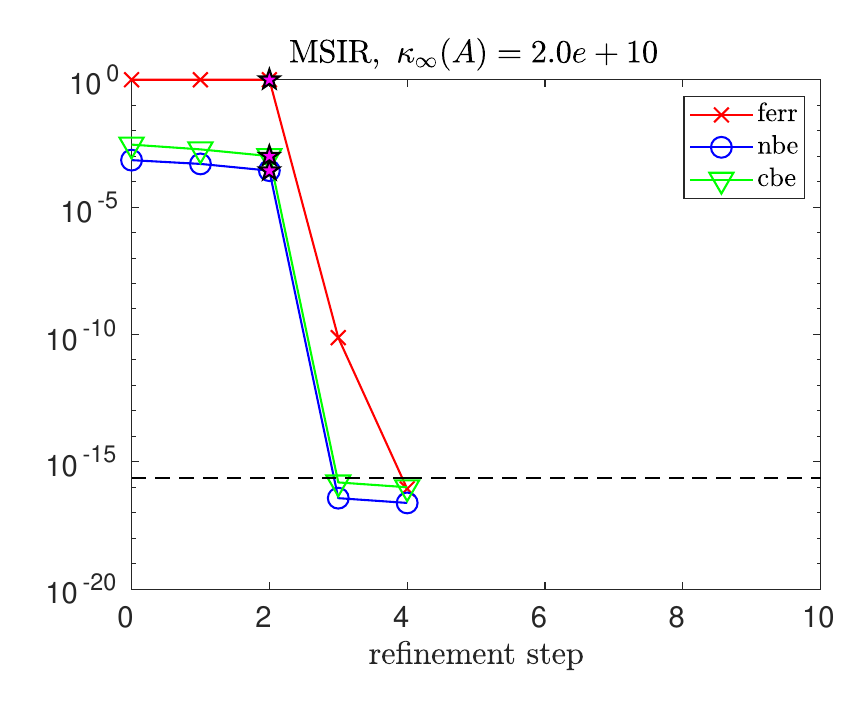}\\
	\includegraphics[trim={3cm 8cm 4cm 8cm},clip, width=.45\textwidth]{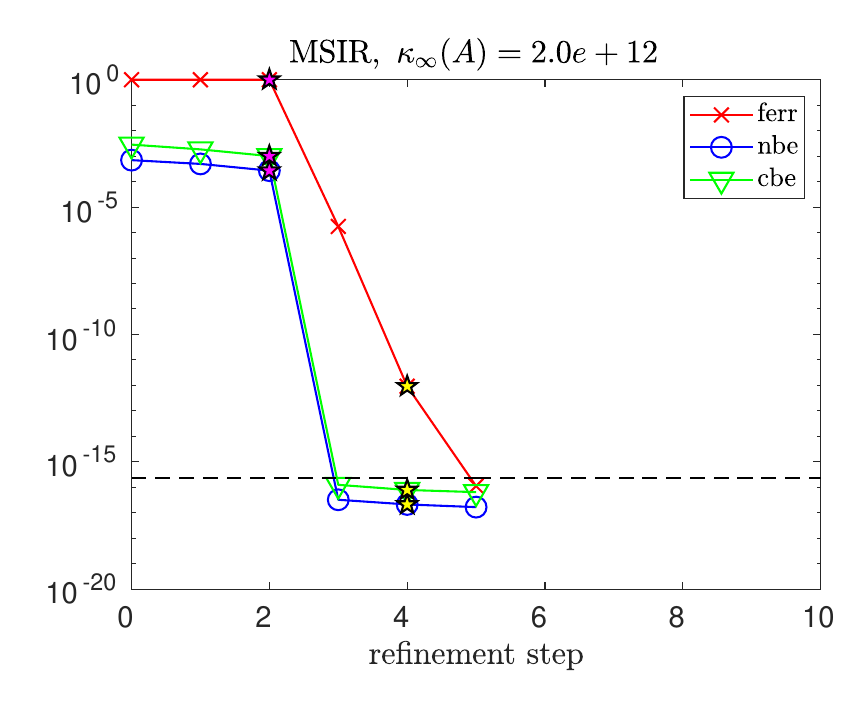}
	\includegraphics[trim={3cm 8cm 4cm 8cm},clip, width=.45\textwidth]{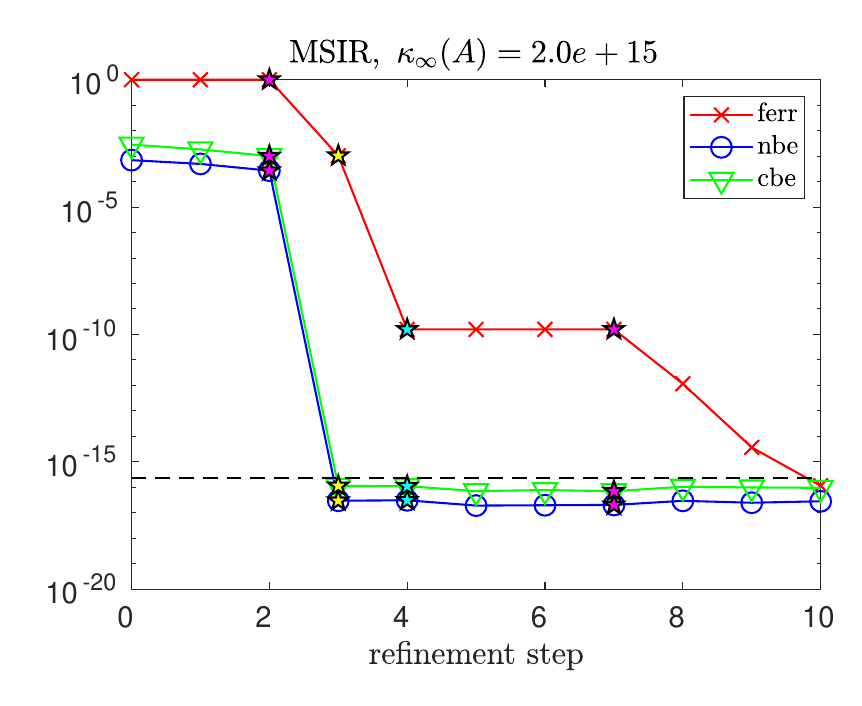}
	\caption{Convergence of errors in MSIR for random dense matrices (mode 2) with $\kappa_2 (A)=10^1$ (top left), $\kappa_2 (A)=10^9$ (top right), $\kappa_2 (A)=10^{11}$ (bottom left), and $\kappa_2 (A)=10^{14}$ (bottom right) for initial precisions $(u_f,u,u_r)$ = (half, double, quad); see also Table \ref{tab:rand_hdq}.}
	\label{fig:randnmat_024}
\end{figure}

Finally, in Table \ref{tab:rand_hdq}, we test initial precisions $(u_f,u,u_r)$ = (half, double, quad). Except for SIR, which again fails to converge once $\kappa_\infty(A)>2\cdot 10^5$, all algorithms converge for $\kappa_\infty(A) < 10^{16}$. Even in cases where SIR converges, the SIR convergence is slow enough that MSIR always switches to SGMRES-IR. 
MSIR only performs the second switch to GMRES-IR for the two most ill-conditioned problems (although this switch is actually unnecessary at least for $\kappa_2(A)=10^{11}$; see also the convergence trajectory in the bottom left plot in Figure \ref{fig:randnmat_024}). Again, this emphasizes the need for a multistage, adaptive approach; if we were to go by the theoretical bounds, we would not expect SGMRES-IR to work beyond $\kappa_\infty(A) = 3\cdot 10^7$, but here it works well even for a matrix with a condition number five orders of magnitude greater. For the most ill-conditioned case, MSIR also performed a precision switch from $(u_f, u, u_r)$ = (half, double, quad) to $(u_f, u, u_r)$ = (single, double, quad) since GMRES-IR required too many GMRES iterations. We plot convergence trajectories for MSIR for the problems with $\kappa_2(A)\in \{10^1, 10^9, 10^{11}, 10^{14}\}$ in Figure \ref{fig:randnmat_024}. 

\begin{table}[h!]
	\centering
	\caption{Number of SIR, SGMRES-IR, GMRES-IR, and MSIR steps with the number of GMRES iterations for each SGMRES-IR and GMRES-IR step for random dense matrices (mode 2) with various condition numbers $\kappa_\infty(A)$, $\kappa_2(A)$, using initial precisions $(u_f,u,u_r)$ = (half, double, quad).}
	\begin{tabular}{|cc|cccc|}
		\hline
		 $\kappa_\infty (A)$     & $\kappa_2(A)$ & {SIR} & {SGMRES-IR} & {GMRES-IR} & {MSIR} \\ \hline
		$2\cdot 10^2$           & $10^1$           & 7                    & (5,5)                      & (5,5)                     & 2, (5)                 \\
		$2\cdot 10^3$           & $10^2$           & 9                    & (5,5)                    & (5,5)                 & 2, (5)              \\
		$2\cdot 10^5$           & $10^4$           & 47                    & (5,6)                    & (5,6)                   & 2, (5,6)               \\
		$2\cdot 10^6$           & $10^5$           & -                    & (5,6)                    & (5,6)                  & 2, (5,6)              \\
		$2\cdot 10^8$           & $10^7$           & -                    & (6,7)                    & (6,7)                   & 2, (6,7)              \\
		$2\cdot 10^{10}$        & $10^9$           & -                    & (7,8)                    & (7,8)                  & 2, (7,8)            \\
		$2\cdot 10^{12}$        & $10^{11}$        & -                    & (8,8,9)                          & (8,9)            & 2, (8,8), (9)          \\
		$2\cdot 10^{15}$        & $10^{14}$        & -                    & (100,100,100,25,10,10)                          & (100,32,10)                         & 2, (10), (10); 3, (4,4,4)                     \\ \hline
	\end{tabular}
	\label{tab:rand_hdq}
\end{table}

\subsubsection{Random dense matrices with geometrically distributed singular values} \label{sec:randnmat_mode3}

We now test the convergence behavior of the IR variants for matrices generated via the MATLAB command \verb|gallery('randsvd',n,kappa(i),3)|. As mentioned, we expect that these are more challenging problems for the GMRES-based iterative refinement schemes since they contain clusters of eigenvalues close to the origin which may fail to be effectively shifted away from the origin by a low precision preconditioner. 

\begin{figure}[h!]
	\centering
	\includegraphics[trim={3cm 8cm 4cm 8cm},clip, width=.45\textwidth]{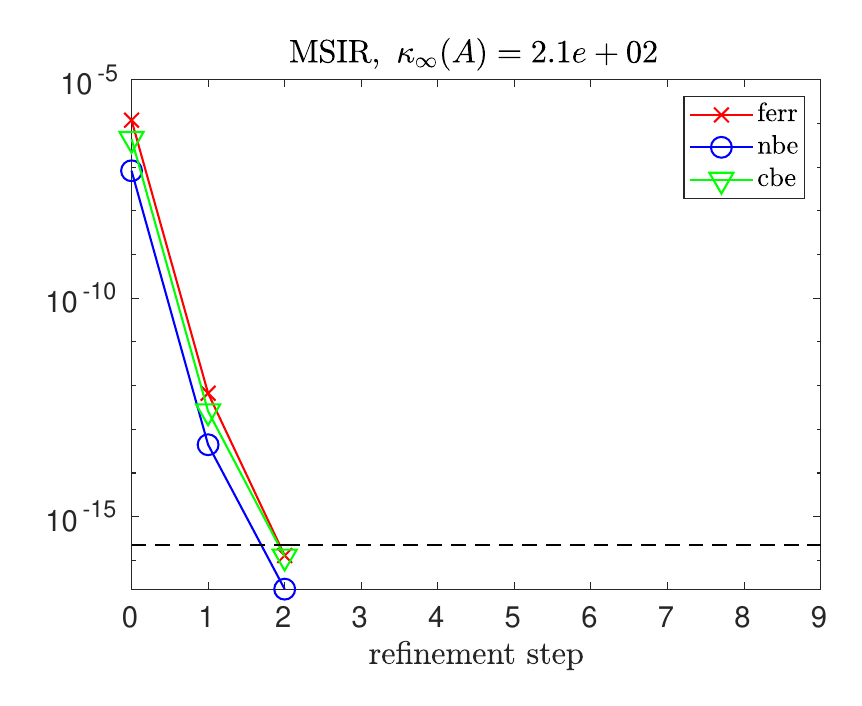}
	\includegraphics[trim={3cm 8cm 4cm 8cm},clip, width=.45\textwidth]{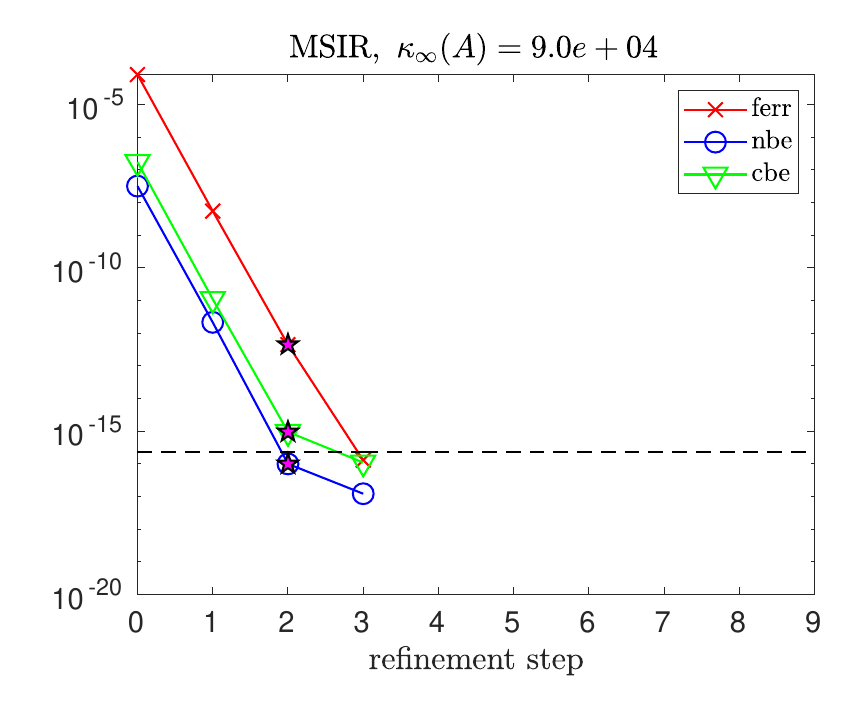}\\
	\includegraphics[trim={3cm 8cm 4cm 8cm},clip, width=.45\textwidth]{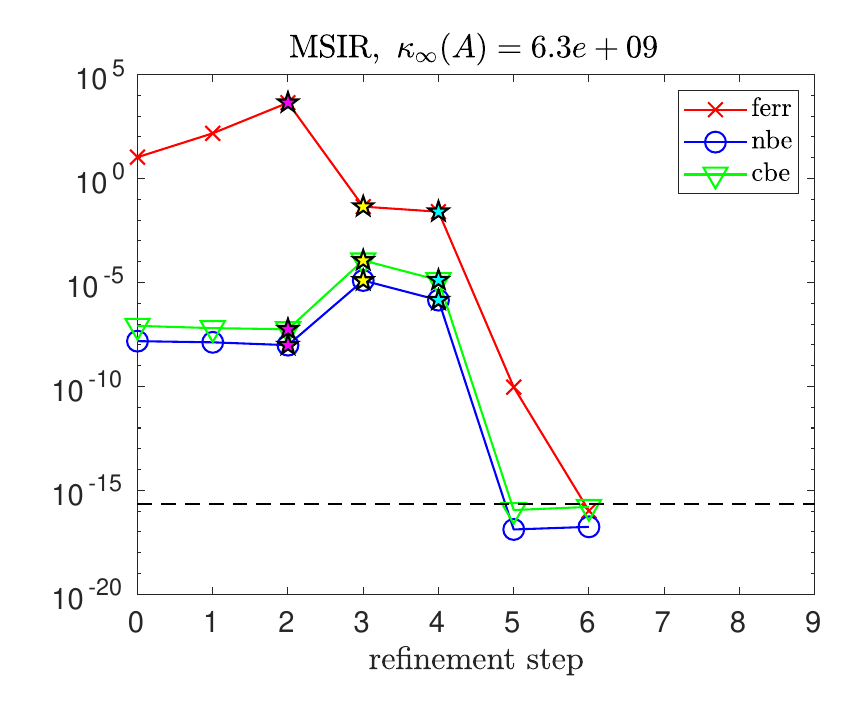}
	\includegraphics[trim={3cm 8cm 4cm 8cm},clip, width=.45\textwidth]{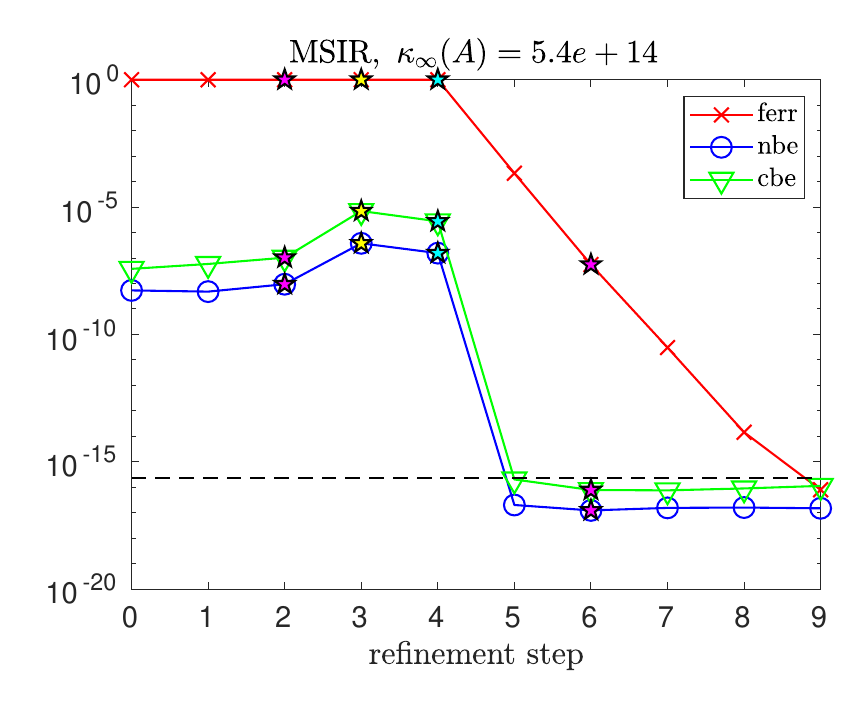}
	\caption{Convergence of errors in MSIR for random dense matrices having geometrically distributed singular values (mode 3) with $\kappa_2 (A)=10^1$ (top left), $\kappa_2 (A)=10^4$ (top right), $\kappa_2 (A)=10^9$ (bottom left), and $\kappa_2 (A)=10^{14}$ (bottom right) for initial precisions $(u_f,u,u_r)$ = (single, double, quad); see also Table \ref{tab:rand_mode3_sdq}.}
	\label{fig:randnmat_mode3_124}
\end{figure}

In Table \ref{tab:rand_mode3_sdq} we test initial precisions $(u_f,u,u_r)$ = (single, double, quad). Here, SIR fails to converge {once $\kappa_\infty(A)$ exceeds $10^9$}. Once $\kappa_\infty(A)$ exceeds $10^9$, the number of GMRES iterations required per SGMRES-IR and GMRES-IR step begins to increase dramatically. This is expected since there is a cluster of eigenvalues remaining close to the origin even after low precision preconditioning. Notice that in these cases, MSIR performs the second switch to GMRES-IR after just one SGMRES-IR step and then performs a precision switch from $(u_f, u, u_r)$ = (single, double, quad) to $(u_f, u, u_r)$ = (double, double, quad) after just one GMRES-IR step, since the number of GMRES iterations per refinement step for both SGMRES-IR and GMRES-IR exceeds our specified value of $k_{max}=0.1n$. We show plots for $\kappa_2(A)\in\{10^1, 10^4, 10^9, 10^{14}\}$ in Figure \ref{fig:randnmat_mode3_124}.

\begin{table}[h!]
	\centering
	\caption{Number of SIR, SGMRES-IR, GMRES-IR, and MSIR steps with the number of GMRES iterations for each SGMRES-IR and GMRES-IR step for random dense matrices having geometrically distributed singular values (mode 3) with various condition numbers $\kappa_\infty(A)$, $\kappa_2(A)$, using initial precisions $(u_f,u,u_r)$ = (single, double, quad).}
	\begin{tabular}{|cc|cccc|}
		\hline
		$\kappa_\infty(A)$     & $\kappa_2(A)$  & {SIR} & {SGMRES-IR} & {GMRES-IR} & {MSIR} \\ \hline
		$2\cdot 10^2$          & $10^1$            & 2                    & (2)                      & (2)                     & 2                     \\
		$10^3$                 & $10^2$            & 2                    & (2)                        & (2)                       & 2                     \\
		$9\cdot 10^4$          & $10^4$            & 3                    & (3)                      & (3)                     & 2, (3)                 \\
		$8\cdot 10^5$          & $10^5$            & 4                    & (4)                      & (4)                     & 2, (4)                 \\
		$7\cdot 10^7$          & $10^7$            & 13                  & (7,7)                    & (7,7)                   & 2, (6,7)               \\
		$6\cdot 10^9$          & $10^9$            & -                    & (22,23,25)                 & (22,26)                & 2, (10), (10); 2        \\
		$6\cdot 10^{11}$       & $10^{11}$         & -                    & (53,54,55)                 & (53,55)                & 2, (10), (10); 2, (2)       \\
		$5\cdot 10^{14}$       & $10^{14}$         & -                    & (84,84,84,85,85)        & (84,88)                & 2, (10), (10); 2, (3,3,3)       \\ \hline
	\end{tabular}
	\label{tab:rand_mode3_sdq}
\end{table}

\begin{figure}[h!]
	\centering
	\includegraphics[trim={3cm 8cm 4cm 8cm},clip, width=.45\textwidth]{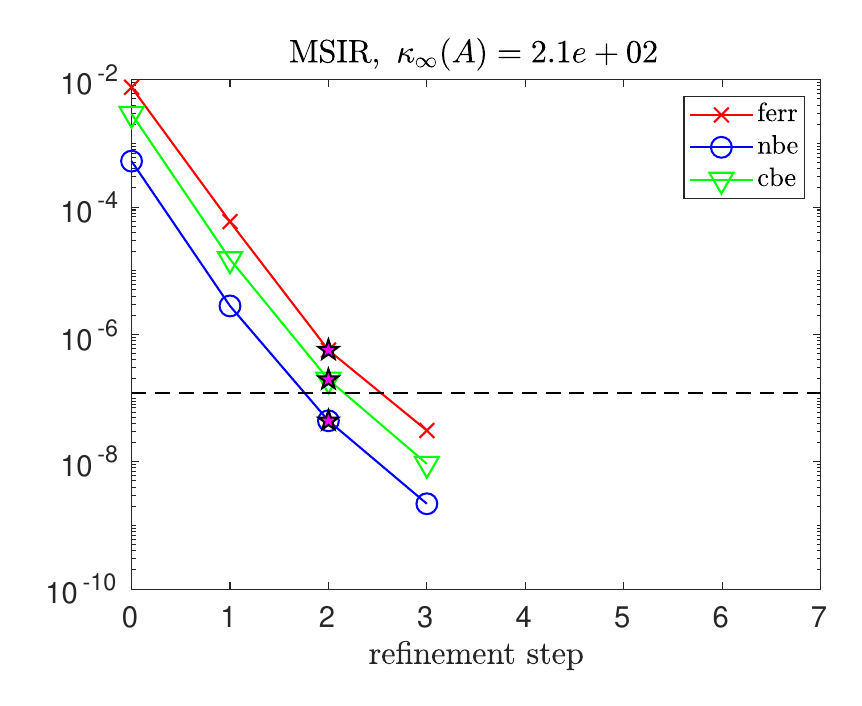}
	\includegraphics[trim={3cm 8cm 4cm 8cm},clip, width=.45\textwidth]{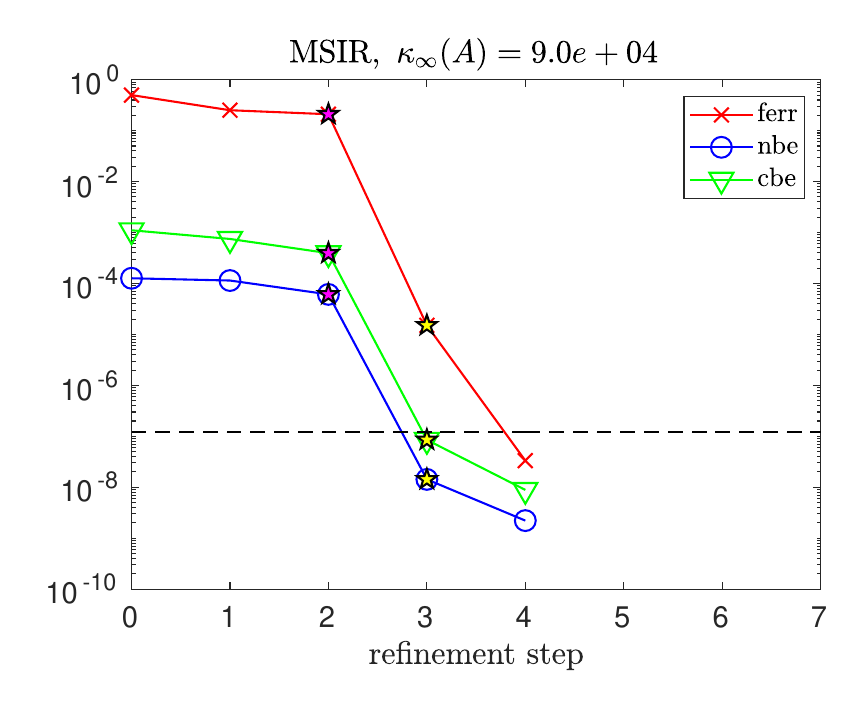}\\
	\includegraphics[trim={3cm 8cm 4cm 8cm},clip, width=.45\textwidth]{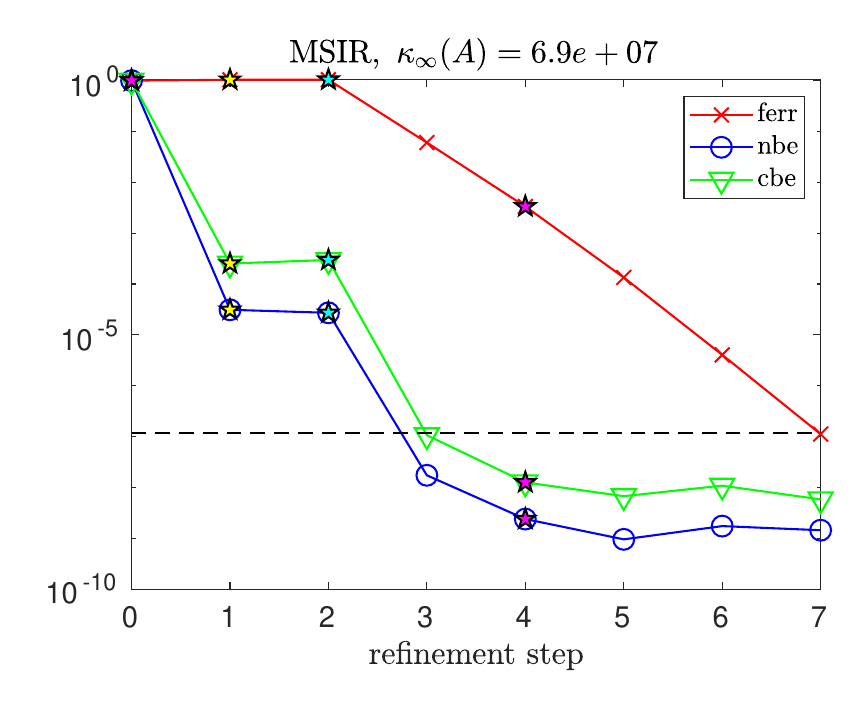}
	\includegraphics[trim={3cm 8cm 4cm 8cm},clip, width=.45\textwidth]{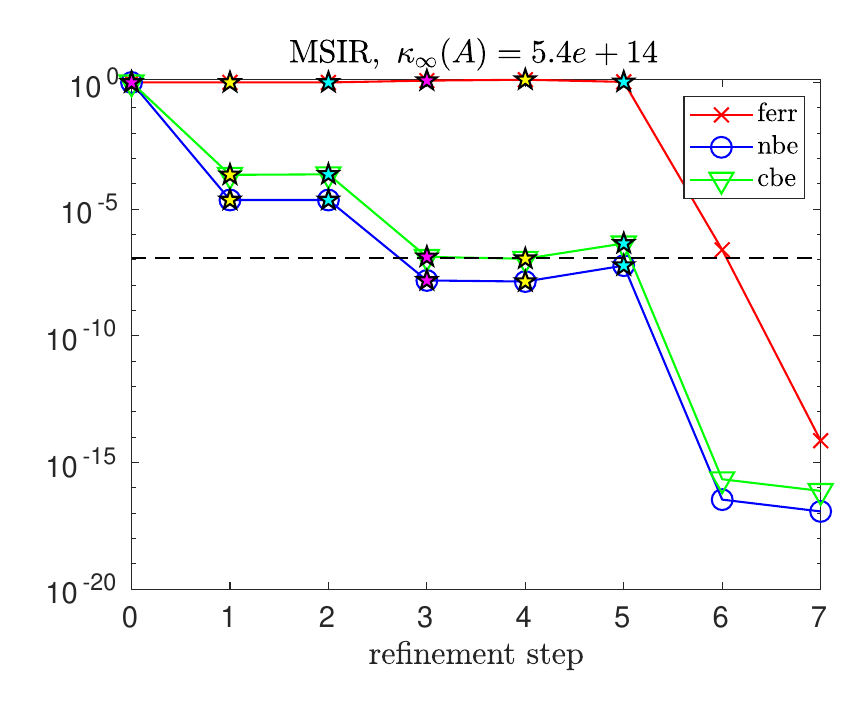}
	\caption{Convergence of errors in MSIR for random dense matrices having geometrically distributed singular values (mode 3) with $\kappa_2 (A)=10^1$ (top left), $\kappa_2 (A)=10^4$ (top right), $\kappa_2 (A)=10^7$ (bottom left), and $\kappa_2 (A)=10^{14}$ (bottom right) for initial precisions $(u_f,u,u_r)$ = (half, single, double); see also Table \ref{tab:rand_mode3_hsd}. }
	\label{fig:randnmat_mode3_012}
\end{figure}

The story for initial precisions $(u_f,u,u_r)$ = (half, single, double), shown in Table \ref{tab:rand_mode3_hsd}, is similar. Once the condition number exceeds around $\sqrt{u}$, the number of GMRES iterations required per refinement step increases significantly for both SGMRES-IR and GMRES-IR. Again, this means that MSIR switches from SGMRES-IR to GMRES-IR and then performs a precision switch (to $(u_f, u, u_r)$ = (single, single, double)) since $k_{max}$ is exceeded. When $\kappa_\infty$ exceeds $10^9$, however, using single precision for the factorization is not enough to combat the slow GMRES convergence, and thus a second precision switch occurs, switching to $(u_f, u, u_r)$ = (double, double, quad). For the most ill-conditioned problem for which SGMRES-IR converges, MSIR does significantly fewer GMRES iterations in total than SGMRES-IR (although we use higher precision for $u_f$ (single), the computational cost will still less considering the total number of GMRES steps performed in SGMRES-IR). Notice also that for $\kappa_2(A)\geq 10^5$, MSIR performs no SIR steps at the initial precision setting $(u_f,u,u_r)$ = (half, single, double) since an Inf or NaN is detected in the first correction term. We plot MSIR convergence for $\kappa_2(A)\in\{10^1, 10^4, 10^7, 10^{14}\}$ in Figure \ref{fig:randnmat_mode3_012}.

\begin{table}[h!]
	\centering
	\caption{Number of SIR, SGMRES-IR, GMRES-IR, and MSIR steps with the number of GMRES iterations for each SGMRES-IR and GMRES-IR step for random dense matrices having geometrically distributed singular values (mode 3) with various condition numbers $\kappa_\infty(A)$, $\kappa_2(A)$, using initial precisions $(u_f,u,u_r)$ = (half, single, double). For $\kappa_2(A) = 10^{11}$, GMRES-IR used (100,100,100,100,100,100,100,100,100,100,99,95,100,97) iterations.}
	\scalebox{0.80}{\begin{tabular}{|cc|cccc|}
		\hline
		$\kappa_\infty(A)$     & $\kappa_2(A)$ & {SIR} & {SGMRES-IR} & {GMRES-IR} & {MSIR} \\ \hline
		$2\cdot 10^2$          & $10^1$            & 3                    & (3)                        & (3)                       & 2, (3)                     \\
		$10^3$                 & $10^2$            & 4                    & (4)                        & (4)                       & 2, (4)                 \\
		$9\cdot 10^4$          & $10^4$            & {\color{black}856}                    & (13,14)                    & (13,14)                   & {\color{black}2, (10), (10)}              \\
		$8\cdot 10^5$          & $10^5$            & -                    & (38,40)                          & (38,41)                         & {\color{black}0, (10), (10); 3}        \\
		$7\cdot 10^7$          & $10^7$            & -                    & (100,100,100,100,100,83)                          & (100,100)                         & {\color{black}0, (10), (10); 2, (5,5,5)}        \\
		$6\cdot 10^9$          & $10^9$            & -                    & -                          & (100,100,100,100)                         & {\color{black}0, (10), (10);1, (10), (10); 1  }                     \\
		$6\cdot 10^{11}$       & $10^{11}$         & -                    & -                          & {\color{black}(100,$\ldots$) \scriptsize{see caption.}}                         & {\color{black}0, (10), (10); 1, (10), (10); 1  }                     \\
		$5\cdot 10^{14}$       & $10^{14}$         & -                    & -                          & -                         & {\color{black}0, (10), (10); 1, (10), (10); 2}                     \\ \hline
	\end{tabular}}
	\label{tab:rand_mode3_hsd}
\end{table}

In Table \ref{tab:rand_mode3_hdq} we show results for initial precisions $(u_f,u,u_r)$ = (half, double, quad). Here the increase in the number of GMRES iterations per GMRES-based refinement step as $\kappa(A)$ increases is also dramatic. One interesting thing to notice here and even moreso in the previous examples with mode 3 matrices is that, in cases where both SGMRES-IR and GMRES-IR converge, the extra precision used in GMRES-IR does not help improve the convergence rate of GMRES (nearly the same number of GMRES iterations per refinement step are required). The aspect that improves is the number of refinement steps. This is likely related to the GMRES convergence trajectories for this class of problems; the residual will nearly stagnate (or at least not make much progress) for a number of iterations and then convergence will happen suddenly. The difference is likely that when this convergence does happen, the extra precision in GMRES-IR results in the approximate solution being found to greater accuracy than in the uniform precision iterations in SGMRES-IR. 

\begin{figure}[h!]
	\centering
	\includegraphics[trim={3cm 8cm 4cm 8cm},clip, width=.45\textwidth]{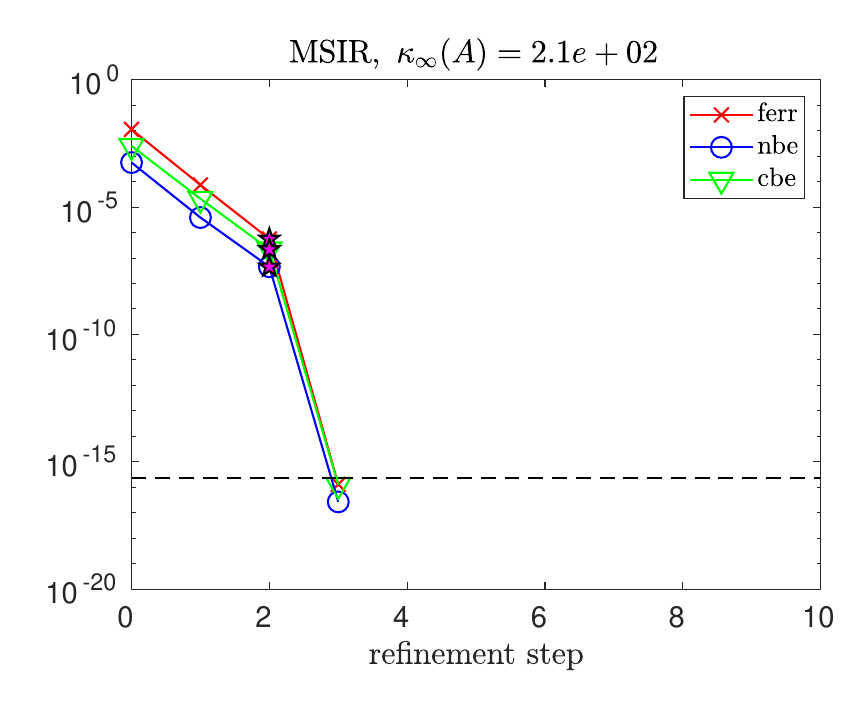}
	\includegraphics[trim={3cm 8cm 4cm 8cm},clip, width=.45\textwidth]{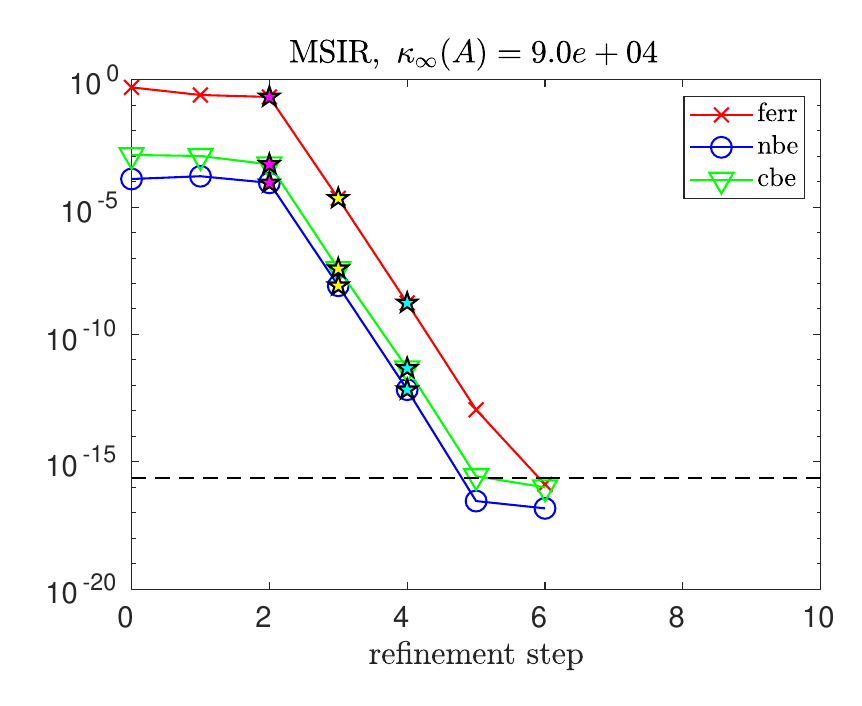}\\
	\includegraphics[trim={3cm 8cm 4cm 8cm},clip, width=.45\textwidth]{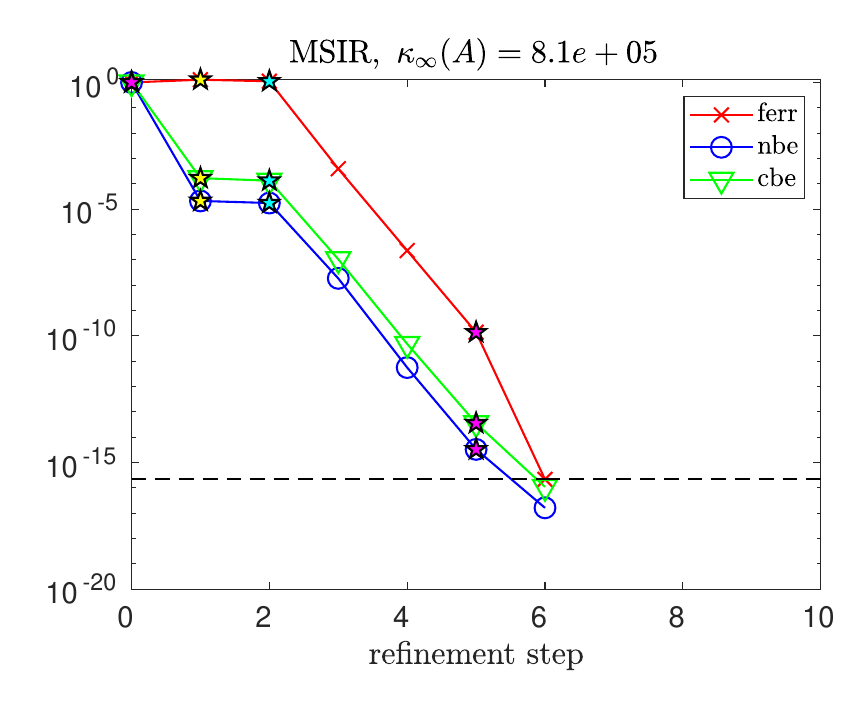}
	\includegraphics[trim={3cm 8cm 4cm 8cm},clip, width=.45\textwidth]{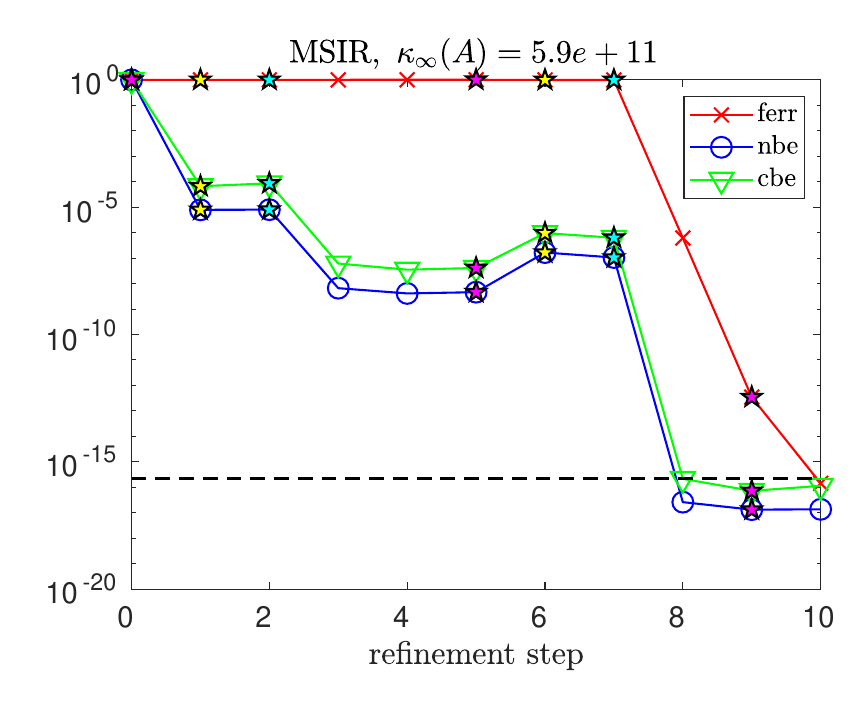}
	\caption{Convergence of errors in MSIR for random dense matrices having geometrically distributed singular values (mode 3) with $\kappa_2 (A)=10^1$ (top left), $\kappa_2 (A)=10^4$ (top right), $\kappa_2 (A)=10^5$ (bottom left), and $\kappa_2 (A)=10^{11}$ (bottom right) for initial precisions $(u_f,u,u_r)$ = (half, double, quad). }
	\label{fig:randnmat_mode3_024}
\end{figure}

{\color{black}As a result of this slow convergence of GMRES, MSIR does a precision switch from $(u_f,u,u_r)$ = (half, double, quad) to $(u_f,u,u_r)$ = (single, double, quad) for matrices with condition number $\kappa_2(A) \geq 10^4$. For matrices with condition number $\kappa_2(A) \geq 10^9$, MSIR does a precision switch a second time, from $(u_f,u,u_r)$ = (single, double, quad) to $(u_f,u,u_r)$ = (double, double, quad). Even though this means that MSIR computes 3 LU factorizations, one in half, one in single, and one in double, this is still likely to be much faster than SGMRES-IR, which does a total of 1841 $O(n^2)$ triangular solves in double precision, or GMRES-IR, which does a total of 927 $O(n^2)$ triangular solves in quadruple precision.  }
We show MSIR convergence for {\color{black}$\kappa_2(A)\in\{10^1, 10^4, 10^5, 10^{11}\}$} in Figure \ref{fig:randnmat_mode3_024}. 

\begin{table}[h!]
	\centering
	\caption{Number of SIR, SGMRES-IR, GMRES-IR, and MSIR steps with the number of GMRES iterations for each SGMRES-IR and GMRES-IR step for random dense matrices having geometrically distributed singular values (mode 3) with various condition numbers $\kappa_\infty(A)$, $\kappa_2(A)$, using initial precisions $(u_f,u,u_r)$ = (half, double, quad). {\color{black}For $\kappa_2(A) = 10^{14}$, SGMRES-IR used (100,100,100,100,100,85,96,93,90,84,94,92,90,81,95,84,89,82,90,96) iterations.}}
	\scalebox{0.8}{\begin{tabular}{|cc|cccc|}
		\hline
		$\kappa_\infty(A)$     & $\kappa_2(A)$  & {SIR} & {SGMRES-IR} & {GMRES-IR} & {MSIR} \\ \hline
		$2\cdot 10^2$          & $10^1$            & 7                    & (5,5)                    & (5,5)                   & 2, (5)               \\
		$10^3$                 & $10^2$            & 9                    & (6,6)                    & (6,6)                   & 2, (6,6)               \\
		$9\cdot 10^4$          & $10^4$            & -                    & (19,20)                 & (19,20)                & {\color{black}2, (10), (10); 2}        \\
		$8\cdot 10^5$          & $10^5$            & -                    & (43,46)                          & (43,46)                         & {\color{black}0, (10), (10); 3, (3)}        \\
		$7\cdot 10^7$          & $10^7$            & -                    & (83,85)                          & (83,85)                         & {\color{black}0, (10), (10); 2, (6,7)}        \\
		$6\cdot 10^9$          & $10^9$            & -                    & (100,100)                          & (100,100)                         & {\color{black}0, (10), (10); 2, (10), (10); 2}                     \\
		$6\cdot 10^{11}$       & $10^{11}$         & -                    & (100,100,100)                          & (100,100)                         & {\color{black}0, (10), (10); 3, (10), (10); 2, (2)}                     \\
		$5\cdot 10^{14}$       & $10^{14}$         & -                    & {\color{black}(100,$\ldots$)}                          & {\color{black}(100,100,100,91,84,88,83,96,85,100)}                         & {\color{black}0, (10), (10); 2, (10), (10); 2, (3,3,3)}                     \\ \hline
	\end{tabular}}
	\label{tab:rand_mode3_hdq}
\end{table}

\subsection{SuiteSparse Matrices} \label{sec:realmat}
\begin{figure}[h!]
	\centering
	\includegraphics[trim={3cm 8cm 4cm 8cm},clip, width=.45\textwidth]{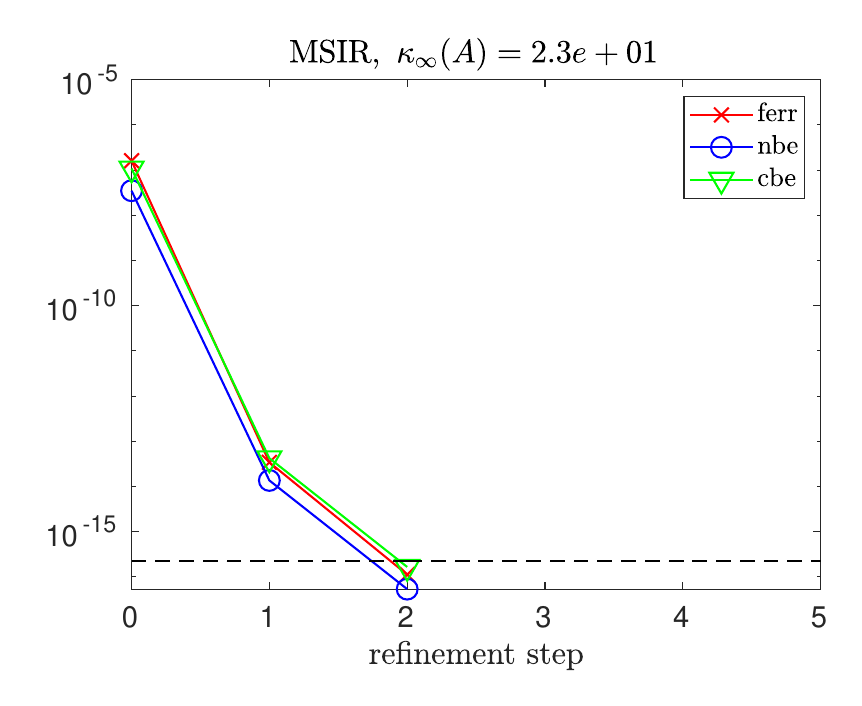}
	\includegraphics[trim={3cm 8cm 4cm 8cm},clip, width=.45\textwidth]{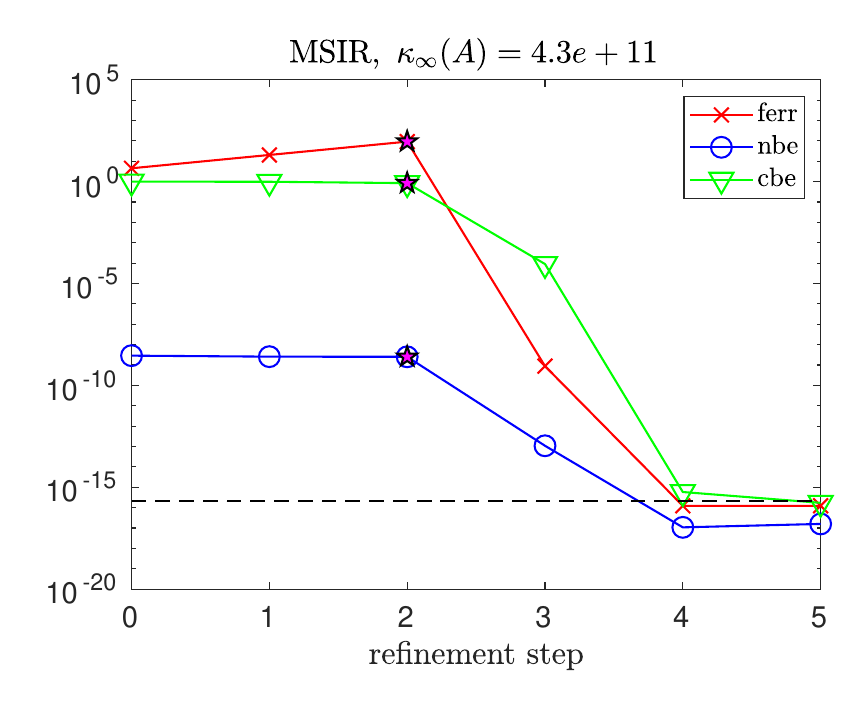}
	\caption{Convergence of errors in MSIR for \textit{cage6} (left) and \textit{ww\_36\_pmec\_36} (right) using initial precisions $(u_f,u,u_r)$ = (single, double, quad).}
	\label{fig:realmat_124}
\end{figure}

\begin{figure}[h!]
	\centering
	\includegraphics[trim={3cm 8cm 4cm 8cm},clip, width=.45\textwidth]{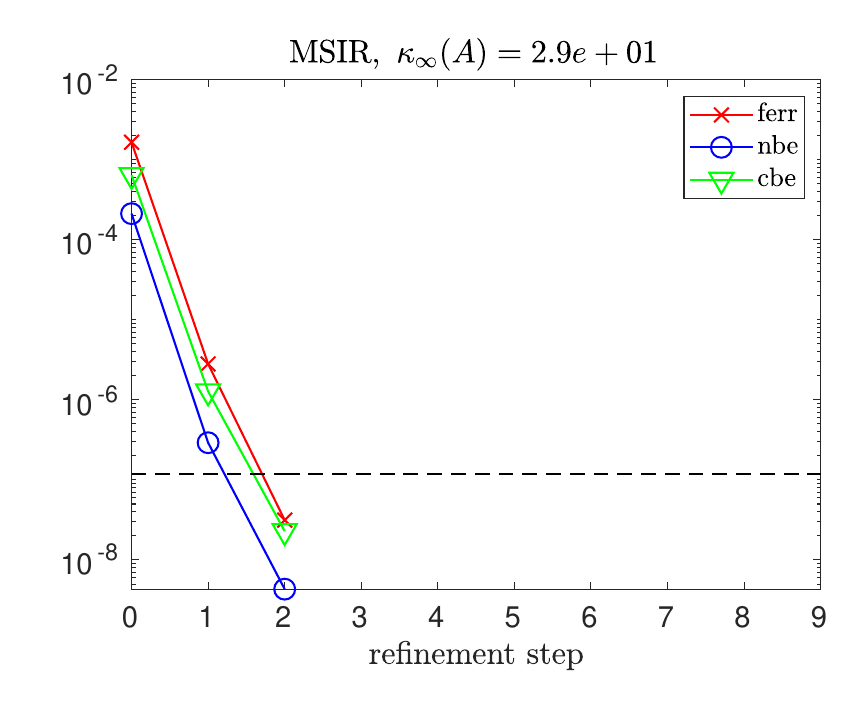}
	\includegraphics[trim={3cm 8cm 4cm 8cm},clip, width=.45\textwidth]{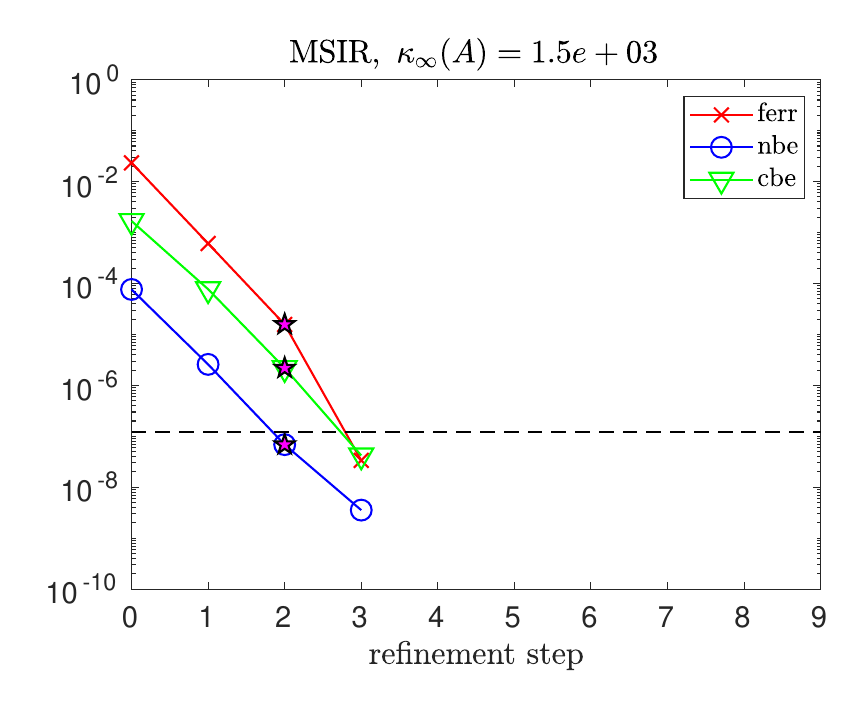}\\
	\includegraphics[trim={3cm 8cm 4cm 8cm},clip, width=.45\textwidth]{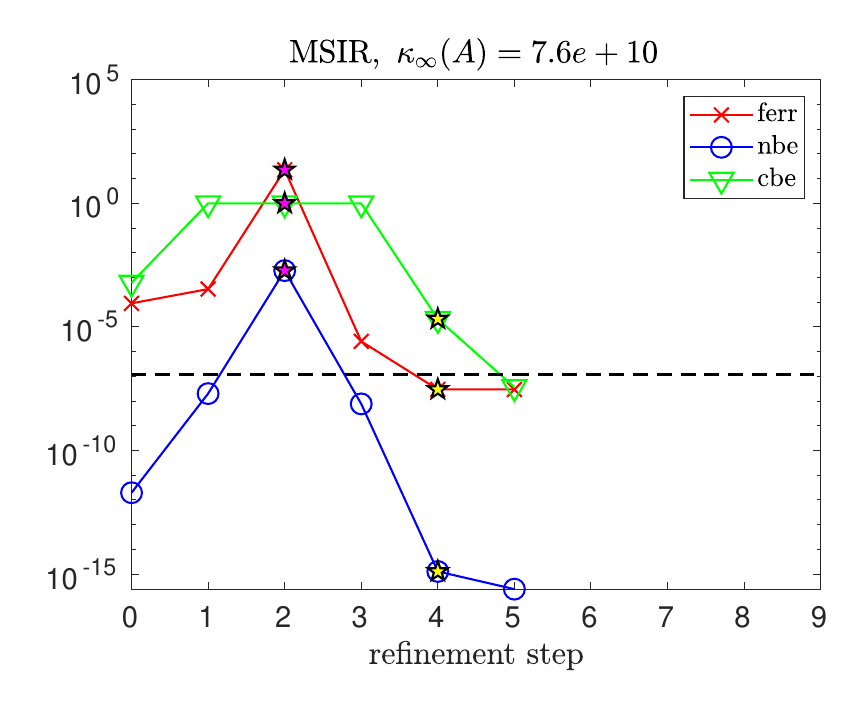}
	\includegraphics[trim={3cm 8cm 4cm 8cm},clip, width=.45\textwidth]{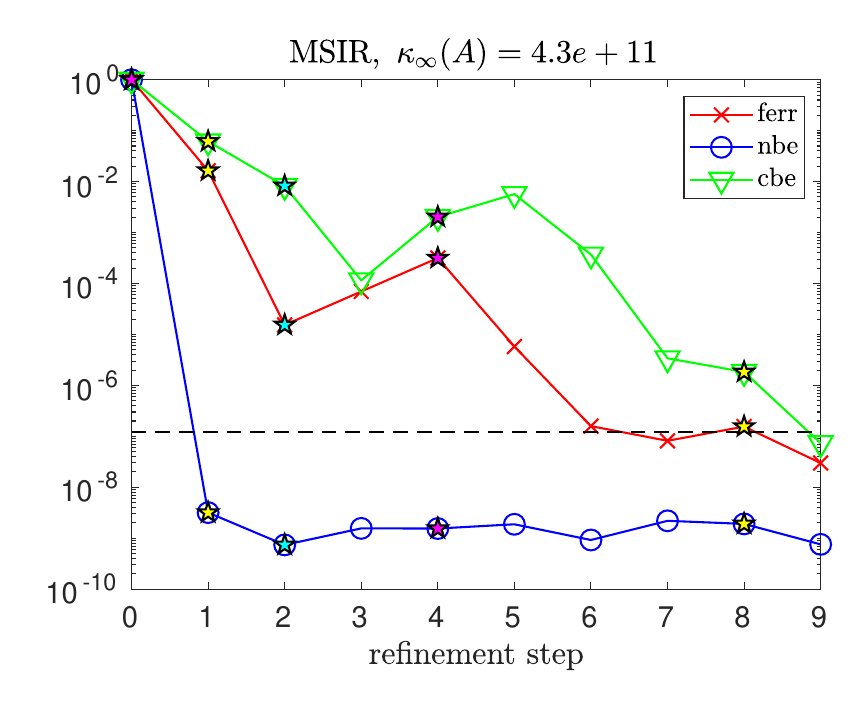}
	\caption{Convergence of errors in MSIR for \textit{cage5} (top left), \textit{bwfa62} (top right), \textit{steam3} (bottom left), and \textit{ww\_36\_pmec\_36} (bottom right) using initial precisions $(u_f,u,u_r)$ = (half, single, double).}
	\label{fig:realmat_012}
\end{figure}

\begin{figure}[h!]
	\centering
	\includegraphics[trim={3cm 8cm 4cm 8cm},clip, width=.45\textwidth]{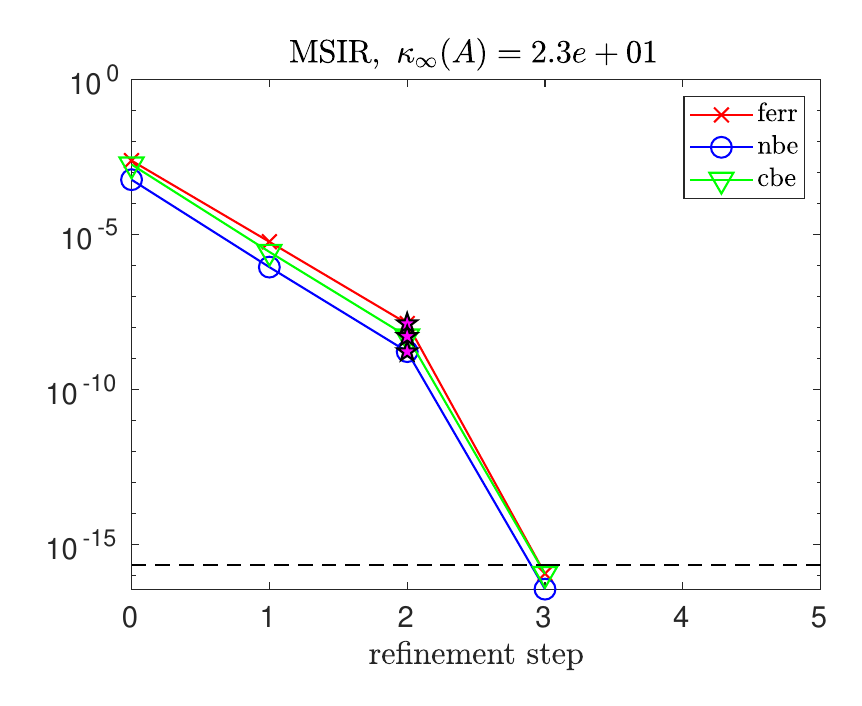}
	\includegraphics[trim={3cm 8cm 4cm 8cm},clip, width=.45\textwidth]{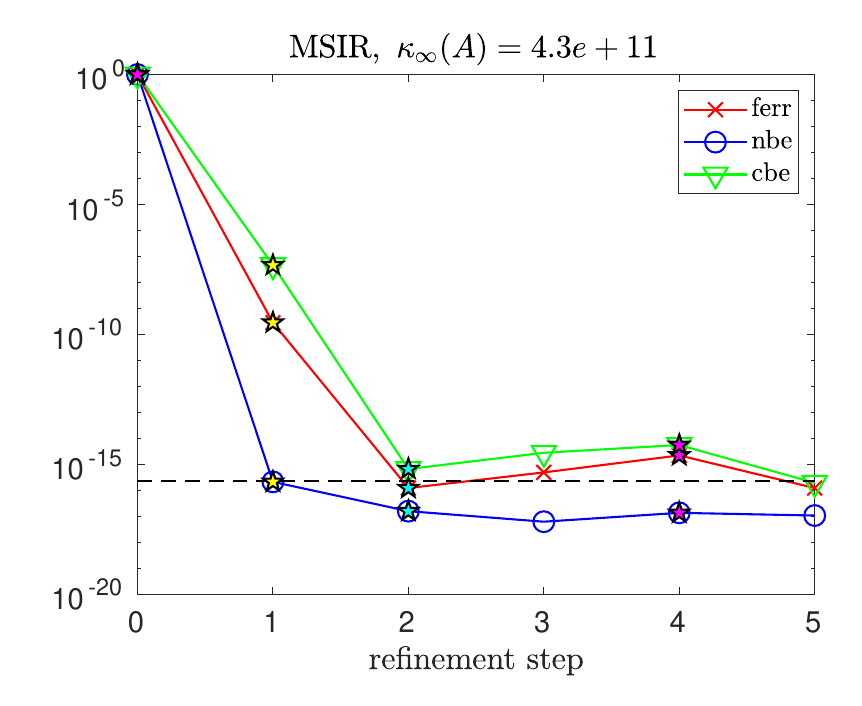}
	\caption{Convergence of errors in MSIR for \textit{cage6} (left) and \textit{ww\_36\_pmec\_36} (right) using initial precisions $(u_f,u,u_r)$ = (half, double, quad).}	
	\label{fig:realmat_024}
\end{figure}

We now test a subset of matrices taken from the SuiteSparse Collection \cite{dh:11} whose properties are listed in Table \ref{tab:matrices}. 
As before, we test SIR, SGMRES-IR, GMRES-IR, and MSIR with {\color{black}initial} precisions $(u_f,u,u_r)$ = (single, double, quad), $(u_f,u,u_r)$ = (half, single, double), and $(u_f,u,u_r)$ = (half, double, quad). For these problems, we always set the right-hand side $b$ to be the vector of ones. 
Tables \ref{tab:real_sdq}-\ref{tab:real_hdq} show the convergence behavior of the different algorithmic variants. For each precision combination, we show a few interesting convergence trajectories in Figures \ref{fig:realmat_124}-\ref{fig:realmat_024}. 

In Table \ref{tab:real_sdq} for initial precisions $(u_f,u,u_r)$ = (single, double, quad), it is observed that SIR converges for all matrices except ww\_36\_pmec\_36  (the most ill-conditioned one). Since SIR converges in few steps, MSIR does not switch to SGMRES-IR. For ww\_36\_pmec\_36, however, MSIR detects nonconvergence of SIR and switches to SGMRES-IR, which then converges in 3 steps. 

\begin{table}[h!]
	\caption{Matrices from \cite{dh:11} used for numerical experiments and their properties.}
	\centering
	\vspace{0.3cm}
	\scalebox{0.89}{
	\begin{tabular}{|c|c|c|c|c|c|}
		\hline
		Name      & Size & Nnz  & $\kappa_\infty$     & Group  & Kind                                \\ \hline
		cage6    & 93   & 785  & 2.34E+01 & vanHeukelum  & Directed Weighted Graph \\
		tols90    & 90   & 1746  & 3.14E+04 & Bai  & Computational Fluid Dynamics Problem \\
		bfwa62    & 62   & 450  & 1.54E+03 & Bai  & Electromagnetics Problem \\
		cage5    & 37   & 233  & 2.91E+01 & vanHeukelum  & Directed Weighted Graph \\
		d\_dyn    & 87   & 230  & 8.71E+06 & Grund  & Chemical Process Simulation Problem \\ 
		d\_ss     & 53   & 144  & 6.02E+08 & Grund  & Chemical Process Simulation Problem \\ 
		Hamrle1   & 32   & 98   & 5.51E+05 & Hamrle & Circuit Simulation Problem          \\ 
		ww\_36\_pmec\_36     & 66  & 1194  & 4.283E+11 & Rommes    & Eigenvalue/Model Reduction Problem    \\
		steam3     & 80  & 314  & 7.64E+10 & HB    & Computational Fluid Dynamics Problem	\\ \hline
	\end{tabular}}
	\label{tab:matrices}
\end{table}

\begin{table}[h!]
	\centering
	\caption{Number of SIR, SGMRES-IR, GMRES-IR, and MSIR steps with the number of GMRES iterations for each SGMRES-IR and GMRES-IR step for real matrices with initial precisions $(u_f,u,u_r)$ = (single, double, quad).}
	\begin{tabular}{|c|cccc|}
		\hline
		Matrix   & SIR & SGMRES-IR & GMRES-IR & MSIR      \\ \hline
		cage6            & 2   & (2)     & (2)    & 2         \\
		tols90           & 2   & (2)     & (2)    & 2         \\
		bfwa62           & 2   & (2)     & (2)    & 2         \\
		cage5            & 2   & (2)     & (2)    & 2         \\
		Hamrle1          & 2   & (2)     & (2)    & 2         \\
		d\_ss            & 2   & (2)       & (2)      & 2         \\
		d\_dyn           & 2   & (2)     & (2)    & 2         \\
		ww\_36\_pmec\_36 & -   & (2,3,3)   & (2,3,3)  & 2, (2,3,3) \\ 
		steam3 			 & 2   & (2)     & (2)    & 2 \\ \hline
	\end{tabular}
	\label{tab:real_sdq}
\end{table}

In Table \ref{tab:real_hsd} for initial precisions $(u_f,u,u_r)$ = (half, single, double), it is seen that SIR converges for the matrices with $\kappa_\infty(A)<10^{6}$. It is also observed that SGMRES-IR fails to converge for the extremely ill-conditioned matrix ww\_36\_pmec\_36. MSIR switches to GMRES-IR for the two most ill-conditioned matrices, although for the steam3 case, the second switch is not technically necessary. It is not clear why for the steam3 matrix, MSIR performs two SGMRES-IR steps and one GMRES-IR step when SGMRES-IR alone only requires one step. {\color{black}For ww\_36\_pmec\_36, MSIR also switches from $(u_f, u, u_r)$ = (half, single, double) to $(u_f, u, u_r)$ = (single, single, double) since GMRES-IR required too many GMRES iterations. It is lastly seen that for steam3, scaling is performed due to the resulting L and U factors containing Inf or NaN when the factorization is performed in half precision, as explained in Section \ref{sec:scaling}.}

\begin{table}[h!]
	\centering
	\caption{Number of SIR, SGMRES-IR, GMRES-IR, and MSIR steps with the number of GMRES iterations for each SGMRES-IR and GMRES-IR step for real matrices with initial precisions $(u_f,u,u_r)$ = (half, single, double). {\color{black}Cases where scaling was required prior to LU factorization are marked with a *..}}
	\begin{tabular}{|c|cccc|}
		\hline
		Matrix   & SIR & SGMRES-IR & GMRES-IR  & MSIR  \\ \hline
		cage6            & 2   & (3)       & (3)       & 2     \\
		tols90           & 2   & (2)       & (2)       & 2     \\
		bfwa62           & 4   & (3)       & (3)       & 2, (3) \\
		cage5            & 2   & (3)       & (3)       & 2     \\
		Hamrle1          & 2   & (2)       & (2)       & 2     \\
		d\_ss            & -   & (3,3)     & (3,3)       & 0, (3,3) \\
		d\_dyn           & -   & (2,2)       & (2,2)       & 0, (2,2) \\
		ww\_36\_pmec\_36 & -   & -         & (66,66,7) & {\color{black}0, (7), (7); 2, (2,2,3,3), (2)}     \\ 
		steam3{\color{black}*}			 & -   & (2)   & (2)  & 2, (2,2), (2) \\ \hline
	\end{tabular}
	\label{tab:real_hsd}
\end{table}

In Table \ref{tab:real_hdq}, using initial precisions $(u_f,u,u_r)$ = (half, double, quad), we see that MSIR is able to converge for all test problems and switches at least to SGMRES-IR in every case due to slow convergence of SIR. For the matrix ww\_36\_pmec\_36, MSIR makes the second switch to GMRES-IR {\color{black}and then does a precision switch from $(u_f, u, u_r)$ = (half, double, quad) to $(u_f, u, u_r)$ = (single, double, quad)} after one step due to the constraint on the number of GMRES iterations per step. {\color{black}As in the $(u_f,u,u_r)$ = (half, single, double) case, due to the use of half precision in the factorization process, scaling is required for steam3.}

Overall, we note that for these real-world sparse problems, to an even greater extent than in the dense test cases, SGMRES-IR and GMRES-IR often still perform well for matrices with condition numbers greatly exceeding constraints given in Table \ref{tab:kappalimit}. We again stress that this motivates our multistage approach, since the theoretical analysis cannot necessarily give indication of which algorithm variant will be the best choice.

\begin{table}[h!]
	\centering
	\caption{Number of SIR, SGMRES-IR, GMRES-IR, and MSIR steps with the number of GMRES iterations for each SGMRES-IR and GMRES-IR step for real matrices with initial precisions $(u_f,u,u_r)$ = (half, double, quad). {\color{black}Cases where scaling was required prior to LU factorization are marked with a *.}}
	\begin{tabular}{|c|cccc|}
		\hline
		Matrix  & SIR & SGMRES-IR & GMRES-IR & MSIR      \\ \hline
		cage6            & 5   & (4,4)     & (4,4)    & 2, (4)     \\
		tols90           & 5   & (3,3)     & (3,3)    & 2, (3)     \\
		bfwa62           & 9   & (4,5)   & (4,5)  & 3, (4)   \\
		cage5            & 5   & (4,4)   & (4,4)  & 2, (3)     \\
		Hamrle1          & 5   & (3,3)     & (3,3)    & 2, (3)     \\
		d\_ss            & -  & (4,5)     & (4,5)    & 0, (4,5)   \\
		d\_dyn           & -   & (4,4)   & (4,4)  & 0, (4,4)   \\
		ww\_36\_pmec\_36 & -   & (8,8)  & (8,9)  & {\color{black}0, (7), (7); 2, (3)} \\ 
		steam3{\color{black}*} 			 & -   & (3,3)   & (3,3)  & 2, (3,3,3) \\ \hline
	\end{tabular}
	\label{tab:real_hdq}
\end{table}

\section{Conclusions and future work}
\label{sec:conclusions}

Mixed precision iterative refinement, and in particular, GMRES-based iterative refinement, {\color{black} has seen great use recently due to the emergence of mixed precision hardware}. The freedom to choose different precision combinations as well as a particular solver leads to an explosion of potential iterative refinement variants, each of which has different performance costs and different condition number constraints under which convergence is guaranteed. {\color{black}In practice, particular iterative refinement variants often work well for problems beyond what is indicated by the analysis, making it difficult for a user to select a priori the least expensive variant that will converge.} 
Motivated by improving usability and reliability, we have developed a multistage mixed precision iterative refinement approach, called MSIR. {\color{black}For a given combination of precisions, MSIR switches between three different iterative refinement variants distinguished by the way in which they solve for the correction to the approximate solution in order of increasing cost, according to stopping criteria adapted from \cite{dh:06}. Then, only if necessary, it increases the precision(s), refactorizes the matrix in a higher precision, and begins again.} We discuss details of the algorithm and perform extensive numerical experiments on both random dense matrices and matrices from SuiteSparse \cite{dh:11} using a variety of {\color{black}initial} precision combinations. Our experiments confirm that the algorithmic variants often outperform what is dictated by the theoretical condition number constraints and demonstrate the benefit of the multistage approach; in particular, our experiments show that there can be an advantage to first trying other solvers before resorting to increasing the precision and refactorizing.   

There are a number of potential extensions to the work presented here. First, while we have only tested IEEE precisions, it is also worthwhile to perform experiments using non-IEEE floating point formats such as bfloat16 \cite{bfloat16}. We also plan to extend the multistage approach to least squares problems, which can be accomplished by combining the error bounds derived in \cite{dh:09} and the three-precision iterative refinement schemes for least squares problems in \cite{ch:20}. Finally, while we have roughly used the number of triangular solves in a given precision as a metric for discussing probable performance characteristics, this may not be a good indication of performance in practice. Thorough high-performance experiments on modern GPUs are needed to appropriately gauge the overhead of the MSIR approach and the algorithm's suitability in practice. 

\bibliographystyle{siamplain}
\bibliography{paper}

\end{document}